\documentclass[a4paper, BCOR=0.0mm, DIV=calc]{scrartcl}
\usepackage[T1]{fontenc}
\usepackage[utf8]{inputenc}
\usepackage[ngerman]{babel}
\usepackage{booktabs, tabularx}
\usepackage[osf,sc]{mathpazo}
\usepackage{textcomp}
\usepackage{microtype}
\usepackage{amsmath, amssymb, amsfonts}
\KOMAoptions{twoside=false, twocolumn=false, headinclude=false, footinclude=false, mpinclude=false, pagesize=auto}
\recalctypearea
\newcommand{\ebinom}[2]{\left(\frac{#1}{#2} \right)}
\begin{document}
\title{\boldmath Analytische Untersuchungen über die Entwicklung der Trinomialpotenz $(1+x+xx)^n$\unboldmath \footnote{
Originaltitel: "`Disquitiones analyticae super evolutione potestatis trinomialis $(1+x+xx)^n$"', erstmals publiziert in "`\textit{Nova Acta Academiae Scientarum Imperialis Petropolitinae} 14, 1805, pp. 75-110"', Nachdruck in "`\textit{Opera Omnia}: Series 1, Volume 16, pp. 56 - 103"', Eneström-Nummer E722, übersetzt von: Alexander Aycock, Textsatz: Artur Diener,  im Rahmen des Projektes "`Eulerkreis Mainz"' }}
\author{Leonhard Euler}
\date{}
\maketitle
\paragraph{§1}
Nachdem ich einst in "`Novorum Commentariorum Tomo XI"' unter dem Titel "`Analytische Beobachtungen"' diese Trinomialpotenz mit viel Eifer untersucht hatte, bin ich auf so außerordentliche Eigenschaften gestoßen, die größerer Aufmerksamkeit der Mathematiker nicht unwürdig schienen.

Deswegen habe ich neulich unternommen denselben Beweis erneut zu behandeln und es haben sich mir, der ich einige Kunstgriffe benutzt habe, viele ausgezeichnete Phänomene offenbart, deren Erklärung ich vertraue den Mathematikern nicht unpassend sein zu werden.
\paragraph{§2}
Ich beginne also mit der Entwicklung selbst dieser Formel
\[
	(1 + x + xx)^n,
\]
die für einzelne Werte des Exponenten $n$ die folgenden Ausdrücke liefert, die in der beigefügten Tabelle dargestellt wurden: \\[1em]
\begin{tabularx}{\textwidth}{lX} 
\toprule
$n$ & $(1+x+xx)^n$ \\ 
\midrule
$0$ & $1$ \\
$1$ & $1 + x + xx$ \\
$2$ & $1 + 2x + 3xx + 2x^3 + x^4$ \\
$3$ & $1 + 3x + 6xx + 7x^3 + 6x^4 + 3x^5 + x^6$ \\
$4$ & $1 + 4x + 10xx + 16x^3 + 19x^4 + 16x^5 + 10x^6 + 4x^7 + x^8$ \\ 
$5$ & $1 + 5x + 15xx + 30x^3 + 45x^4 + 51x^5 + 45x^6 + 30x^7 + 15x^8 + 5x^9 + x^{10}$ \\
$\mathrm{etc.}$ & $\mathrm{etc.}$ \\
\bottomrule
\end{tabularx}
Hier wird natürlich aus einer beliebigen Potenz sehr leicht das folgende geschlossen; wenn nämlich für einen beliebigen Wert des Exponenten $n$ ein beliebiger Koeffizient mit den beiden vorhergehenden in eine Summe gezogen wird, erhält man den Koeffizienten für die folgende unten hinzuschreibende Potenz des Exponenten $n+1$.
\paragraph{§3}
Dem diese Tabelle anzuschaunden ist sofort klar, dass in einer beliebigen Entwicklung die Koeffizienten der Terme bis hin zur Mitte, die Größe $x^n$ meint, wachsen, daher aber wiederum in der selben Ordnung schrumpfen bis hin zum letzten Term, der $x^{2n}$ ist. Darauf durchschaut man auch nicht schwerm dass für die Potenz $(1 + x + xx)^n$ im Allgemeinen die Anfangsterme so ausgedrückt werden können:
\begin{align*}
	&1 + nx + \frac{n(n+1)}{1 \cdot 2}x^2 + \frac{n(n-1)(n+4)}{1 \cdot 2 \cdot 3}x^3 + \frac{n(n-1)(nn + 7n - 6)}{1 \cdot 2 \cdot 3 \cdot 4}x^4 \\
	&\frac{n(n-1)(n-2)(n+12)}{1 \cdot 2 \cdot 3 \cdot 4 \cdot 5}x^5 + \mathrm{etc.}
\end{align*}
Es kommt aber darauf nicht an, dass diese Terme weiter fortschreiten, weil in deren Koeffizienten keine Struktur entdeckt wird.
\paragraph{§4}
Hier aber achte ich besonders auf den größten oder den mittleren Koeffizienten, den ich für die Potenz $(1+x+xx)^n$ im Allgemeinen immer gleich $px^n$ setzen werde; dann aber werde ich die daher folgenden Terme so darstellen: $qx^{n+1},\, rx^{n+2},\, sx^{n+3},\, tx^{n+4},\, \mathrm{etc}$, woher die vorhergehenden Mittelterme, indem man der Reihe nach zurückgeht, $qx^{n-1},\, rx^{n-2},\, sx^{n-3},\, tx^{n-4},\, \mathrm{etc.}$ sein werden.  Darauf aber werde ich für die folgende Potenz $(1+x+xx)^{n+1}$ die selben Buchstaben mit einem Strich versehen, natürlich $p',\, q',\, r',\, s',\, \mathrm{etc}$, welche ich weiter für die erneut folgende Potenz $(1+x+xx)^{n+2}$ mit einem Doppelstrich bezeichnen werde; für die folgenden mit einem dreifachen Strich, einem vierfachen und so weiter.
\paragraph{§5}
Nachdem diese Dinge vorausgeschickt worden sind, werde ich in dieser Abhandlung aus den oberen Reihen der Tabelle hauptsächlich die Mittelterme, die mit den größten Koeffizienten versehen worden sind, betrachten, welche $1,\, x,\, 3x^2,\, 7x^3,\, 19x^4,\, 51x^5,\, \mathrm{etc.}$ sind, die zusammengenommen eine Reihe festlegen, deren Summe ich mit dem Buchstaben $P$ kennzeichnen werde, sodass 
\[
	P = 1 + x + 3x^2 + 7x^3 + 19x^4 + 51x^5 + \cdots + px^n + p'x^{n+1} + p''x^{n+2} + \mathrm{etc.}
\]
ist.
\paragraph{§6}
Außerdem aber, so wie diese Terme aus der oberen Tabelle entlang der Diagonale ausgesucht worden sind, wollen wir auf ähnliche Weise eine solche Reihe entlang der weiteren jener parallelen Diagonalen bilden, die Summen welcher Reihen wir in gleicher Weise mit speziellen Buchstaben auf die folgende Weise bezeichnen wollen:
\begin{align*}
Q &= x^2 &+ 2x^3 &+ 6x^4 &+ 16x^5 &+ 45x^6 &+ \cdots &+ qx^{n+1} &+ q'x^{n+2} &+ q''x^{n+3} &+ \mathrm{etc.} \\
R &= x^4 &+ 3x^5 &+ 10x^6 &+ 30x^7 &+ \cdots & &+ rx^{n+2} &+ r'x^{n+3} &+ r''x^{n+4} &+ \mathrm{etc.} \\
S &= x^6 &+ 4x^7 &+ 15x^8 &+ \cdots & & &+ sx^{n+3} &+ s'x^{n+4} &+ s''x^{n+5} &+ \mathrm{etc.} \\
T &= x^8 &+ 5x^9 &+ \cdots & & & &+ tx^{n+4} &+ t'x^{n+5} &+ t''x^{n+6} &+ \mathrm{etc.} \\
\mathrm{etc.} 
\end{align*}
Nachdem das alles festgesetzt wurde, ist mir zuerst vorgelegt, die Werte der Kleinbuchstaben $p,\, q,\, r,\, s,\, \mathrm{etc.}$ und deren Derivate $p',\, q',\, r',\, s',\, \mathrm{etc.}\, p'',\, q'',\, r'',\, s'',\, \mathrm{etc.}$ zu finden; nachdem das gemacht worden ist, werde ich auch die Werte der Großbuchstaben $P,\, Q,\, R,\, S,\, \mathrm{etc.}$ erforschen.
\paragraph{§7}
Weil $p$ der Koeffizient aus der Entwicklung der Formel $(1+x+xx)^n$ zu entstehenden Potenz $x^n$ ist, wollen wir diese Formel auf diese Weise darstellen:
\[
	(x(1+x) + 1)^n,
\]
für die Entwicklung welcher wir die schon einige Male genutzte Bezeichnungsweise gebrauchen wollen, durch die ich die Koeffizienten der gleichen Binomialpotenz durch diese Charaktere zu bezeichnen pflege, $\ebinom{n}{1},\, \ebinom{n}{2},\, \ebinom{n}{3},\, \ebinom{n}{4},\, \ebinom{n}{5},\, \mathrm{etc}$, sodass \\
\begin{align*} 
\ebinom{n}{1} &= n \\
\ebinom{n}{2} &= \frac{n(n-1)}{1 \cdot 2} \\
\ebinom{n}{3} &= \frac{n(n-1)(n-2)}{1 \cdot 2 \cdot 3} \\
\ebinom{n}{4} &= \frac{n(n-1)(n-2)(n-3)}{1 \cdot 2 \cdot 3 \cdot 4} \\
\ebinom{n}{5} &= \frac{n(n-1)(n-2)(n-3)(n-4)}{1 \cdot 2 \cdot 3 \cdot 4 \cdot 5} \\
\vdots &  \\
\ebinom{n}{\lambda} &= \frac{n(n-1)(n-2)(n-3) \cdots (n - \lambda + 1)}{1 \cdot 2 \cdot 3 \cdots \lambda}
\end{align*}
Über diese Charaktere wird es hier förderlich sein bemerkt zu haben, dass im Allgemeinen immer
\[
	\ebinom{n}{\lambda} = \ebinom{n}{n - \lambda}
\]
ist, weil diese Koeffizienten rückwärts dieselbe Ordnung beibehalten; und weil die äußersten Koeffizienten die Einheit sind, wird
\[
	\ebinom{n}{0} = \ebinom{n}{n} = 1
\]
sein. Darauf, weil aus dem Gesetz der Progression so alle Terme, die dem ersten vorangehen, wie die Terme, die dem letzten folgen, verschwinden, wird es sein wie folgt:
\begin{align*}
\ebinom{n}{-1} &= \ebinom{n}{n+1} = 0 \\
\ebinom{n}{-2} &= \ebinom{n}{n+2} = 0 \\
\ebinom{n}{-3} &= \ebinom{n}{n+3} = 0 \\
\mathrm{etc.}
\end{align*}
\paragraph{§8}
Nachdem diese Dinge vorausgeschickt worden sind, wird unsere Formel $(x(1+x) + 1)^n$ auf gewohnte Weise wie ein Binom entwickelt diese Reihe geben:
\[
	x^n(1+x)^n + \ebinom{n}{1}x^{n-1}(1+x)^{n-1} + \ebinom{n}{2}x^{n-2}(1+x)^{n-2} + \ebinom{n}{3}x^{n-3}(1+x)^{n-3} + \mathrm{etc},
\]
wo man bemerke, dass im Allgemeinen
\[
	(1+x)^{\lambda} = 1 + \ebinom{\lambda}{1}x + \ebinom{\lambda}{2}x^2 + \ebinom{\lambda}{3}x^3 + \mathrm{etc.}
\]
ist. Aus diesen einzelnen Gliedern jener erläuterten Formel müssen die Terme, die die Potenz $x^n$ enthalten, herausgenommen werden, die natürlich zusammengenommen den Mittelterm $px^n$ zusammensetzen.
\paragraph{§9}
Das erste Glied aber, $x^n(1+x)^n$, liefert nur den Term $x^n$ dieser Form. Aus dem zweiten Glied aber wird der zweite Term diese Form haben, die $\ebinom{n}{1}\ebinom{n-1}{1}x^n$ ist. Aus dem dritten Glied entsteht die Potenz $x^n$ aus dem dritten Term, welcher $\ebinom{n}{2}\ebinom{n-2}{2}x^n$ ist. Auf ähnliche Weise folgert man aus dem vierten Glied $\ebinom{n}{3}\ebinom{n-3}{3}x^n$. Aus dem fünften Glied entsteht $\ebinom{n}{4}\ebinom{n-4}{4}x^n$ und so weiter. Daher wird also der wahre Wert des Buchstaben $p$ so berechnet:
\[
	p = 1 + \ebinom{n}{1}\ebinom{n-1}{1} + \ebinom{n}{2}\ebinom{n-2}{2} + \ebinom{n}{3}\ebinom{n-3}{3} + \ebinom{n}{4}\ebinom{n-4}{4} + \mathrm{etc.}
\]
\paragraph{§10}
Auf ähnliche Weise lassen sich aus derselben Entwicklung die Koeffizienten der Potenz $x^{n+1}$ berechnen, die zusammengenommen den Wert des Buchstaben $q$ geben werden. Eine solche Potenz aber, die aus dem ersten Glied entsteht, wird $\ebinom{n}{1}x^{n+1}$ sein. Aus dem zweiten Glied entsteht $\ebinom{n}{1}\ebinom{n-1}{2}x^{n+1}$, aus dem dritten Glied $\ebinom{n}{2}\ebinom{n-2}{3}x^{n+1}$, aus dem vierten $\ebinom{n}{3}\ebinom{n-3}{4}x^{n+1}$ und so weiter, weshalb der wahre Wert des Buchstaben $q$ auf diese Weise ausgedrückt werden wird:
\[
	q = \ebinom{n}{1} + \ebinom{n}{1}\ebinom{n-1}{2} + \ebinom{n}{2}\ebinom{n-2}{3} + \ebinom{n}{3}\ebinom{n-3}{4} + \mathrm{etc.}
\]
wo wegen der Analogie der erste Term, $\ebinom{n}{1}$, so dargestellt zu verstehen ist: $\ebinom{n}{0}\ebinom{n}{1}$. Wenn nämlich, weil ein beliebiger Term aus zwei Faktoren besteht, die ersten Terme diese Reihe festsetzen: $\ebinom{n}{0},\, \ebinom{n}{1},\, \ebinom{n}{2},\, \ebinom{n}{3},\, \ebinom{n}{4},\, \mathrm{etc}$, setzen die hinteren $n$ aber diese fest: $\ebinom{n}{1},\, \ebinom{n-1}{2},\, \ebinom{n-2}{3},\, \ebinom{n-3}{4},\, \mathrm{etc.}$
\paragraph{§11}
Auf die gleiche Weise wird aus den Potenzen $x^{n+2}$, die aus den einzelnen Gliedern gefolgert werden, der Term $rx^{n+2}$ gebildet werden; aber in der Tat liefert das erste Glied für diese Potenz $1 \cdot \ebinom{n}{2}x^{n+2}$ oder der Analogie wegen $\ebinom{n}{0}\ebinom{n}{2}x^{n+2}$. Aus dem zweiten Glied entsteht dieselbe Potenz $\ebinom{n}{1}\ebinom{n-1}{3}x^{n+2}$, aus dem dritten Glied $\ebinom{n}{2}\ebinom{n-2}{4}x^{n+2}$, aus dem vierten $\ebinom{n}{3}\ebinom{n-3}{5}x^{n+2}$ und so weiter; aus diesen gesammelten also erhalten wir den Wert des Buchstaben $r$ auf diese Weise ausgedrückt:
\[
	r = \ebinom{n}{0}\ebinom{n}{2} + \ebinom{n}{1}\ebinom{n-1}{3} + \ebinom{n}{2}\ebinom{n-2}{4} + \ebinom{n}{3}\ebinom{n-3}{5} + \mathrm{etc.}
\]
\paragraph{§12}
Es wäre überflüssig, dieselbe Deduktion die folgenden Buchstaben hinzuzufügen, weil ja schon hinreichend klar ist, dass
\begin{align*}
s &= \ebinom{n}{0}\ebinom{n}{3} + \ebinom{n}{1}\ebinom{n-1}{4} + \ebinom{n}{2}\ebinom{n-2}{5} + \ebinom{n}{3}\ebinom{n-3}{6} + \mathrm{etc.} \\ 
t &= \ebinom{n}{0}\ebinom{n}{4} + \ebinom{n}{1}\ebinom{n-1}{5} + \ebinom{n}{2}\ebinom{n-2}{6} + \ebinom{n}{3}\ebinom{n-3}{7} + \mathrm{etc.} \\
u &= \ebinom{n}{0}\ebinom{n}{5} + \ebinom{n}{1}\ebinom{n-1}{6} + \ebinom{n}{2}\ebinom{n-2}{7} + \ebinom{n}{3}\ebinom{n-3}{8} + \mathrm{etc.} \\
\mathrm{etc.}
\end{align*}
sein wird und es wird im Allgemeinen, wenn wir der Potenz $x^{n + \lambda}$ den Buchstaben $z$ zuteilen,
\[
	z = \ebinom{n}{0}\ebinom{n}{\lambda} + \ebinom{n}{1}\ebinom{n-1}{\lambda + 1} + \ebinom{n}{2}\ebinom{n-2}{\lambda + 2} + \ebinom{n}{3}\ebinom{n-3}{\lambda + 3} + \mathrm{etc.} \\
\]
sein.
\paragraph{§13}
Es ist hier klar, dass alle Terme dieser Reihen in dieser allgemeinen Form $\ebinom{n}{\alpha}\ebinom{n-\alpha}{\beta}$ enthalten sind, welche ich bemerke immer gleich zu sein zu $\ebinom{n}{\beta}\ebinom{n-\beta}{\alpha}$, sodass die Buchstaben $\alpha$ und $\beta$ eine Permutation zulassen. Weil nämlich, nachdem die Entwicklung gemacht worden ist,
\[
	\ebinom{n}{\alpha} = \frac{n(n-1)(n-2)(n-3) \cdots (n-\alpha + 1)}{1 \cdot 2 \cdot 3 \cdots \alpha}
\]
ist und
\[
	\ebinom{n-\alpha}{\beta} = \frac{(n-\alpha)(n-\alpha - 1)(n-\alpha -2) \cdots (n-\alpha -\beta +1)}{1 \cdot 2 \cdot 3 \cdot 4 \cdots \beta}
\]
und nach der Multiplikation wird
\[
	\ebinom{n}{\alpha}\ebinom{n-\alpha}{\beta} = \frac{n(n-1)(n-2)(n-3) \cdots (n-\alpha -\beta +1)}{1\cdot 2\cdot 3\cdots \alpha \cdot 1 \cdot 2 \cdot 4 \cdots \beta}
\]
sein, wo die Vertauschbarkeit der Buchstaben $\alpha$ und $\beta$ ins Auge springt.
\paragraph{§14}
Wenn also die vorher gefundenen Reihen gleich auf diese Weise umgewandelt werden, lässt die erste für $p$ gefundene freilich keine Änderung zu, und die übrigen werden auf folgende Weise ausgedrückt:
\begin{align*}
q &= \ebinom{n}{1}\ebinom{n-1}{0} + \ebinom{n}{2}\ebinom{n-2}{1} + \ebinom{n}{3}\ebinom{n-3}{2} + \ebinom{n}{4}\ebinom{n-4}{3} + \mathrm{etc.} \\
r &= \ebinom{n}{2}\ebinom{n-2}{0} + \ebinom{n}{3}\ebinom{n-3}{1} + \ebinom{n}{4}\ebinom{n-4}{2} + \ebinom{n}{5}\ebinom{n-5}{3} + \mathrm{etc.} \\
s &= \ebinom{n}{3}\ebinom{n-3}{0} + \ebinom{n}{4}\ebinom{n-4}{1} + \ebinom{n}{5}\ebinom{n-5}{2} + \ebinom{n}{6}\ebinom{n-6}{3} + \mathrm{etc.} \\
\vdots \\
z &= \ebinom{n}{\lambda}\ebinom{n-\lambda}{0} + \ebinom{n}{\lambda +1}\ebinom{n-\lambda -1}{1} + \ebinom{n}{\lambda +2}\ebinom{n-\lambda -2}{2} + \mathrm{etc.} \\
\end{align*}
\paragraph{§15}
Außerdem aber ist diese Transformation besonders bemerkenswert, durch die
\[
	\ebinom{n}{\alpha}\ebinom{n-\alpha}{\beta} = \ebinom{\alpha + \beta}{\alpha}\ebinom{n}{\alpha + \beta}
\]
ist. Weil nämlich
\[
	\ebinom{\alpha + \beta}{\alpha} = \frac{(\alpha + \beta)(\alpha + \beta - 1)(\alpha + \beta -2) \cdots (\beta + 1)}{1 \cdot 2 \cdot 3 \cdots \alpha}
\]
ist und
\[
	\ebinom{n}{\alpha + \beta} = \frac{n(n-1)(n-2) \cdots (n-\alpha - \beta + 1)}{1 \cdot 2 \cdot 3 \cdots \alpha \bullet (\alpha +1)(\alpha +2)\cdots (\alpha + \beta)}
\]
oder
\[
	\ebinom{n}{\alpha + \beta} = \frac{n(n-1)(n-2) \cdots (n -\alpha - \beta + 1)}{1 \cdot 2 \cdot 3 \cdots \beta \bullet (\beta + 1) \cdots (\beta + 2)(\beta + \alpha)},
\]
wird das Produkt
\[
\ebinom{\alpha + \beta}{\alpha}\ebinom{n}{\alpha + \beta} = \frac{n(n-1)(n-2) \cdots (n-\alpha - \beta + 1)}{1 \cdot 2 \cdot 3 \cdots \beta \bullet 1 \cdot 2 \cdot 3 \cdots \alpha},
\]
in welche selbe Form die Formel $\ebinom{n}{\alpha}\ebinom{a-\alpha}{\beta}$ aufgelöst wird.
\paragraph{§16}
Durch diese Transformation also können die oberen Reihen auf folgende Weise ausgedrückt werden:
\begin{align*}
p &= 1 + \ebinom{2}{1}\ebinom{n}{2} + \ebinom{4}{2}\ebinom{n}{4} + \ebinom{6}{3}\ebinom{n}{6} + \mathrm{etc.} \\
q &= \ebinom{1}{0}\ebinom{n}{1} + \ebinom{3}{1}\ebinom{n}{3} + \ebinom{5}{2}\ebinom{n}{5} + \ebinom{7}{3}\ebinom{n}{7} + \mathrm{etc.} \\
r &= \ebinom{2}{0}\ebinom{n}{2} + \ebinom{4}{1}\ebinom{n}{4} + \ebinom{6}{2}\ebinom{n}{6} + \ebinom{8}{3}\ebinom{n}{8} + \mathrm{etc.} \\
s &= \ebinom{3}{0}\ebinom{n}{3} + \ebinom{5}{1}\ebinom{n}{5} + \ebinom{7}{2}\ebinom{n}{7} + \ebinom{9}{3}\ebinom{n}{9} + \mathrm{etc.} \\ 
\vdots \\
z &= \ebinom{\lambda}{0}\ebinom{n}{\lambda} + \ebinom{\lambda +2}{1}\ebinom{n}{\lambda +2} + \ebinom{\lambda +4}{2}\ebinom{n}{\lambda +4} + \ebinom{\lambda +6}{3}\ebinom{n}{\lambda +6} + \mathrm{etc.} \\
\end{align*}
\paragraph{§17}
Es verdient noch eine andere Transformation bemerkt zu werden, die bei numerischer Rechnung besonders geeignet ist. Weil nämlich aus der ersten Form
\[
	z = \ebinom{n}{\lambda} + \ebinom{n}{1}\ebinom{n-1}{\lambda + 1} + \ebinom{n}{2}\ebinom{n-2}{\lambda + 2} + \mathrm{etc,}
\]
ist, ist ein beliebiger Term dieser Reihe $\ebinom{n}{\alpha}\ebinom{n-\alpha}{\lambda +\alpha}$, der gleich $\Pi$ genannt werde, und es wird, nachdem die Entwicklung gemacht worden ist,
\[
	\Pi = \frac{n(n-1)(n-2) \cdots (n-\alpha +1) \cdots (n-2\alpha - \lambda +1)}{1 \cdot 2 \cdot 3 \cdots \alpha \bullet 1\cdot 2\cdot 3 \cdots (\lambda + \alpha)}
\]
sein. Wenn wir also anstelle von $\alpha$ hier $\alpha +1$ schreiben, dass der folgende Term entsteht, der also gleich
\[
	\frac{n(n-1)(n-2) \cdots (n-2\alpha - \lambda -1)}{1 \cdot 2 \cdot 3 \cdots (\alpha + 1) \bullet 1 \cdot 2 \cdot 3 \cdots (\lambda +\alpha +1)}
\]
sein wird, liefert dieser durch jenen geteilt den Quotienten
\[
	\frac{(n-2\alpha - \lambda)(n -2\alpha - \lambda -1)}{(\alpha +1)(\lambda + \alpha +1)}.
\]
Daher also wird der folgende Term
\[
	\Pi\frac{(n-2\alpha - \lambda)(n -2\alpha - \lambda -1)}{(\alpha +1)(\lambda + \alpha +1)}
\]
sein.
\paragraph{§18}
Wenn also nach Newton'scher Manier der Buchstabe $\Pi$ einen beliebigen vorhergehenden Term bezeichnet, wird der folgende immer 
\[
	\Pi\frac{(n-2\alpha -\lambda)(n-2\alpha -\lambda -1)}{(\alpha +1)(\lambda + \alpha +1)}
\]
sein; weil daher der erste Term $\ebinom{n}{\lambda}$ ist, wo $\alpha = 0$ ist, wird, wenn dieser durch $\Pi$ bezeichnet wird, der zweite Term gleich
\[
	\Pi\frac{(n-\lambda)(n-\lambda - 1)}{(\lambda +1)}
\]
sein; wenn dieser erneut $\Pi$ genannt wird, wird der dritte Term gleich
\[
	\Pi\frac{(n-\lambda -2)(n-\lambda -3)}{2(\lambda +2)}
\]
sein; wenn dieser erneut $\Pi$ genannt wird, wird der vierte Term gleich
\[
	\Pi\frac{(n-\lambda -4)(n-\lambda -5)}{3(\lambda +3)}
\]
sein und so weiter. Auf diese Weise wird unsere Reihe für $z$ diese Form annehmen:
\begin{align*}
	z = \ebinom{n}{\lambda} &+ \Pi\frac{(n-\lambda)(n-\lambda - 1)}{1(\lambda +1)} + \Pi\frac{(n-\lambda -2)(n-\lambda -3)}{2(\lambda +2)} \\
	&+ \Pi\frac{(n-\lambda -4)(n-\lambda -5)}{3(\lambda +3)} + \mathrm{etc},
\end{align*}
wo natürlich $\Pi$ immer den vorhergehenden Term bezeichnet.
\paragraph{§19}
Daher also, wenn wir anstelle von $\lambda$ nacheinander die Wert $0,\, 1,\, 2,\, 3,\, \mathrm{etc}$ schreiben, werden wir für unsere Buchstaben $p,\, q,\, r,\, s,\, \mathrm{etc}$ die folgenden Reihen erhalten:
\begin{align*}
p &= 1 + \Pi\frac{n(n-1)}{1\cdot 1} + \Pi\frac{(n-2)(n-3)}{2\cdot 2} + \Pi\frac{(n-4)(n-5)}{3 \cdot 3} + \mathrm{etc.} \\
q &= \ebinom{n}{1} + \Pi\frac{(n-1)(n-2)}{1\cdot 2} + \Pi\frac{(n-3)(n-4)}{2\cdot 3} + \Pi\frac{(n-5)(n-6)}{3\cdot 4} + \mathrm{etc.} \\
r &= \ebinom{n}{2} + \Pi\frac{(n-2)(n-3)}{1\cdot 3} + \Pi\frac{(n-4)(n-5)}{2\cdot 4} + \Pi\frac{(n-6)(n-7)}{3\cdot 5} + \mathrm{etc.} \\
s &= \ebinom{n}{3} + \Pi\frac{(n-3)(n-4)}{1\cdot 4} + \Pi\frac{(n-5)(n-6)}{2\cdot 5} + \Pi\frac{(n-7)(n-8)}{3\cdot 6} + \mathrm{etc.} \\
\mathrm{etc.} 
\end{align*}
\paragraph{§20}
Diese Formen sind zur numerischen Berechnung besonders geeignet, was genügen wird, das für den Buchstaben $p$ allein gezeigt zu haben. Wir suchen natürlich eines Beispiels wegen den Wert von $p$ für den Fall $n=6$, und man wird seine einzelnen Teile auf die folgende Weise finden:\\
\begin{align*} 
I. &= 1 = 1 \\
II. &= 1\cdot \frac{6\cdot 5}{1\cdot 1} = 30 \\
III. &= 30\cdot \frac{4\cdot 3}{2\cdot 2} = 90 \\
IV. &= 90\cdot\frac{2\cdot 1}{3\cdot 3} = 20 \\
\end{align*}
also ist die Summe gleich $p=141$.
\paragraph{§21}
Auf ähnliche Weise wollen wir den Wert von $p$ für den Fall $n=12$ suchen, deren einzelne Teile man auf folgende Weise berechnen wird:
\begin{align*}
I.   &= 1 = 1 \\
II.  &= 1     \cdot\frac{12\cdot 11}{1\cdot 1} = 132 \\
III. &= 132   \cdot\frac{10\cdot 9}{2\cdot 2}  = 2970 \\
IV.  &= 2970  \cdot\frac{8\cdot 7}{3\cdot 3}   = 18480 \\
V.   &= 18480 \cdot\frac{6\cdot 5}{4\cdot 4}   = 34650 \\
VI.  &= 34650 \cdot\frac{4\cdot 3}{5\cdot 5}   = 16632 \\
VII. &= 16632 \cdot\frac{2\cdot 1}{6\cdot 6}   = 924 \\
\end{align*}
also ist die Summe gleich $p=73789$.
\paragraph{§22}
Bald aber werden wir eine um Vieles angenehmere Art angeben, die einzelnen Terme dieser Reihen aus den beiden vorhergehenden zu finden, woher mit leichter Rechnung alle Werte für die Buchstaben $p,\, q,\, r,\, \mathrm{etc}$ für die einzelnen Exponenten $n$ beschafft werden können; und so werden sich alle diese Werte, wie weit es beliebt, fortsetzen lassen. Wir wollen aber diese Relation zuerst getrennt für die Zahlen, die unter dem Buchstaben $p$ enthalten sind, aufstellen.
\section*{\boldmath Untersuchung der Relation zwischen den drei aufeinander folgenden Werten $p,\, p',\, p''$ \unboldmath}
\paragraph{§23}
Weil
\[
	p = 1 + \ebinom{n}{1}\ebinom{n-1}{1} + \ebinom{n}{2}\ebinom{n-2}{2} +\ebinom{n}{3}\ebinom{n-3}{3} +\mathrm{etc.}
\]
ist, wollen wir einen gewissen Term $\ebinom{n}{\alpha}\ebinom{n-\alpha}{\alpha}$ dieser Reihe betrachten, den wir gleich $\Pi$ nennen wollen, sodass, nachdem die Entwicklung gemacht worden ist,
\[
	\Pi = \frac{n(n-1)(n-2)\cdots (n-2\alpha +1)}{1 \cdot 2 \cdot 3 \cdots \alpha \bullet 1 \cdot 2 \cdot 3 \cdots \alpha}
\]
ist; den Term aber, der diesem folgt, wollen wir durch $\Phi$ bezeichnen, sodass
\[
	\Phi = \ebinom{n}{\alpha +1}\ebinom{n-\alpha -1}{\alpha +1}
\]
ist und daher nach Ausführung der Entwicklung
\[
	\Phi = \frac{n(n-1)(n-2) \cdots (n-2\alpha -1)}{1\cdot 2\cdot 3 \cdots (\alpha +1) \bullet 1\cdot 2\cdot 3\cdots (\alpha +1)};
\]
und daher wird man also
\[
	\frac{\Phi}{\Pi} = \frac{(n-2\alpha)(n-2\alpha -1)}{(\alpha + 1)(\alpha +1)} \quad \text{und daher} \quad \Pi = \frac{(\alpha +1)(\alpha +1)}{(n-2\alpha)(n-2\alpha -1)}\Phi
\]
haben.

\paragraph{§24}
Gleich wollen wir für die folgenden Werte $p'$ und $p''$ die entsprechenden Werte von $\Phi$ durch $\Phi'$ und $\Phi''$ bezeichnen; weil diese ja aus dem Wert $\Phi$ entstehen, wenn $n+1$ und $n+2$ anstelle von $n$ geschrieben wird, wird nach Ausführung der Entwicklung
\[
	\Phi' = \frac{(n+1)n(n-1) \cdots (n-2\alpha)}{1\cdot 2\cdot 3 \cdots (\alpha + 1) \bullet 1 \cdot 2 \cdot 3 \cdots (\alpha + 1)}
\] 
sein, woher klar ist, dass
\[
	\Phi' : \Phi = \frac{n+1}{n-2\alpha -1}
\]
sein wird und daher
\[
	\Phi' = \frac{n+1}{n-2\alpha -1}\Phi.
\]
Auf ähnliche Weise, wenn wir auch hier $n+1$ anstelle von $n$ schreiben, werden wir
\[
	\Phi'' = \frac{n+2}{n-2\alpha}\Phi' \quad \text{oder} \quad \Phi'' = \frac{(n+1)(n+2)\Phi}{(n-2\alpha-1)(n-2\alpha)}
\]
haben.
\paragraph{§25}
Daher wollen wir gleich diesen Ausdruck bilden:
\[
	A\Phi + \frac{B}{n+1}\Phi' + \frac{C}{(n+1)(n+2)}\Phi'',
\]
deren Wert also durch den Buchstaben $\Phi$ selbst so ausgedrückt werden wird:
\[
	\Phi \left( A + \frac{B}{n-2\alpha -1} + \frac{C}{(n-2\alpha -1)(n-2\alpha)}\right),
\]
wo wir versuchen wollen, die Buchstaben $A,\, B,\, C$ so zu bestimmen, dass diese Form den vorhergehenden Term $\Pi$ selbst gleich wird. Es ist aber klar, dass diese Buchstaben, damit sie auf alle Terme in gleicher Weise erstreckt werden können, den Buchstaben $\alpha$ nicht involvieren müssen. Nachdem also anstelle von $\Pi$ der vor durch $\Phi$ ausgedrückte gegebene Wert eingesetzt wurde, werden wir die folgende durch $\Phi$ geteilte Gleichung haben:
\[
	A + \frac{B}{n-2\alpha -1} + \frac{C}{(n-2\alpha -1)(n-2\alpha)} = \frac{(\alpha +1)(\alpha +1)}{(n-2\alpha -1)(n -2\alpha)},
\]
die von Brüchen befreit
\[
	A(n-2\alpha +1)(n-2\alpha) + B(n-2\alpha) + C = (\alpha +1)(\alpha +1)
\]
wird.
\paragraph{§26}
Weil in dieser Gleichung der Buchstabe $\alpha$ zum zweiten Grad aufsteigt, werden die drei Buchstaben $A,\, B,\, C$ im Allgemeinen ausreichen, dass sie aus dieser Gleichung bestimmt werden können. Zuerst also wollen wir auf beiden Seiten die Terme, die das Quadrat $\alpha\alpha$ involvieren, gleichsetzen, woher diese Gleichung entstehen wird:
\[
	4A\alpha\alpha = \alpha\alpha \quad \text{und daher} \quad A = \tfrac{1}{4}.
\]
Auf dieselbe Weise wollen wir die Terme gleichsetzen, die den Buchstaben $\alpha$ in sich involvieren, woher wir auf diese Gleichung geführt werden:
\[
	2\alpha(1-2n)A - 2\alpha B = 2\alpha,
\]
woher
\[
	B = -\frac{2n-1}{4} - 1 = -\frac{2n-3}{4}
\]
wird. Schließlich geben die von $\alpha$ unabhängigen Terme diese Gleichung:
\[
	(nn-n)A + Bn + C = 1,
\]
woher man
\[
	C = \frac{(n+2)^2}{4}
\]
findet.
\paragraph{§27}
Nachdem also diese Werte gefunden worden sind, wird für die einzelnen Terme immer
\[
	A\Phi + \frac{B}{n+1}\Phi' + \frac{C}{(n+2)(n+1)}\Phi'' = \Pi
\]
sein. Wenn wir daher also daraus diese Formel berechnen:
\[
	Ap + \frac{B}{n+1}p' + \frac{C}{(n+2)(n+1)}p'',
\]
wird aus den ersten für $\Phi$ angenommenen Termen die vorhergehenden der Reihe $p$ entstehen, welche gleich $0$ ist; aus dem zweiten für $\Phi$ angenommenen Term wird aber der erste Term, der $1$ ist; aus den dritten Termen aber wird der zweite Term gebildet, welcher $\ebinom{n}{1}\ebinom{n-1}{1}$ ist; aus den vierten für $\Phi$ angenommenen Termen entsteht der dritte, der $\ebinom{n}{2}\ebinom{n-2}{2}$ ist, und so weiter; und so werden die drei auf diese Weise berechneten Reihen diese Reihe ergeben:
\[
	0 + 1 + \ebinom{n}{1}\ebinom{n-1}{1} + \ebinom{n}{2}\ebinom{n-2}{2} + \mathrm{etc,}
\]
welches die für $p$ gegebene Reihe selbst ist. Daher werden wir zwischen den drei Buchstaben $p,\, p',\, p''$ diese Gleichung haben:
\[
	Ap + \frac{B}{n+1}p' + \frac{C}{(n+2)(n+1)}p'' = p.
\]
\paragraph{§28}
Wir wollen nun anstelle der Buchstaben $A,\, B,\, C$ die gerade gefundenen Werte einsetzen und unsere Gleichung zwischen diesen drei Buchstaben [Formel], die zurückgeführt wird auf diese [Formel], woher [Formel] wird.
\paragraph{§29}
Daher können also leicht für die einzelnen Werte des Exponenten $n$ alle mit dem Buchstaben $p$ bezeichneten Zahlen bestimmt werden, während ein beliebiger aus zwei vorhergehenden zusammengesetzt ist. So wird für $n=0$ genommen $p=1$ und $p'=1$ sein und daher der dritte
\[
	p'' = 1 + \tfrac{1}{2}(1 + 3\cdot 1) = 3.
\]
Für $n=1$ genommen wird $p=1$, $p'=3$ sein und daher der vierte Term 
\[
	p'' = 3 + \tfrac{2}{3}(3 + 3\cdot 1) = 7.
\]
Für $n=2$ genommen wird wegen $p=3$ und $p'=7$ der fünfte Term
\[
	p'' = 7 + \tfrac{3}{4}(7 + 3\cdot 3) = 19
\]
sein. Wenn $n=3$ genommen wird, wird wegen $p=7$ und $p'=19$ der sechste Term
\[
	p'' = 19 + \tfrac{4}{5}(19 + 3\cdot 7) = 51
\]
sein.
\paragraph{§30}
Wenn wir auf diese Weise weiter vorgehen, werden wir diese Progression, wie es beliebt, mithilfe der Form
\[
	p' + \frac{n+1}{n+2}(p' + 3p) = p''
\]
fortsetzen können, die uns folgende Bestimmungen liefert:
\begin{align*}
	51 + \tfrac{5}{6}(51 + 3\cdot 19) &= 141 \\
	141 + \tfrac{6}{7}(141 + 3\cdot 51) &= 393 \\
	393 + \tfrac{7}{8}(393 + 3\cdot 141) &= 1107 \\
	1107 + \tfrac{8}{9}(1107 + 3\cdot 393) &= 3139 \\
	3139 + \tfrac{9}{10}(3139 + 3\cdot 1107) &= 8953 \\
	8953 + \tfrac{10}{11}(8953 + 3\cdot 3139) &= 25653 \\
	25653 + \tfrac{11}{12}(25653 + 3\cdot 8953) &= 73789 \\
	\mathrm{etc.}	
\end{align*}
\paragraph{§31}
Mit der gleichen Methode kann auch die Relation zwischen drei aufeinander folgenden Termen für die folgenden Buchstaben $q,\, r,\, s,\, \mathrm{etc}$ gefunden werden. Damit wir diese Sache allgemeiner machen, wollen wir die Relation zwischen den drei Termen $z,\, z'$ und $z''$ untersuchen, von denen wir annehmen wollen, dass ihnen der Buchstabe $\lambda$ entspricht.
\section*{\boldmath Untersuchung der Relation zwischen drei aufeinander folgenden Termen $z,\, z',\, z''$ \unboldmath}
\paragraph{§32}
Weil
\[
	z = \ebinom{n}{\lambda} + \ebinom{n}{1}\ebinom{n-1}{\lambda +1} + \ebinom{n}{2}\ebinom{n-2}{\lambda + 2} + \ebinom{n}{3}\ebinom{n-3}{\lambda +3} + \mathrm{etc.}
\]
ist, wollen wir den allgemeinen Term dieser Reihe
\[
	\Pi = \ebinom{n}{\alpha}\ebinom{n-\alpha}{\lambda + \alpha}
\]
betrachten, dessen Wert entwickelt
\[
	\Pi = \frac{n(n-1)(n-2) \cdots (n - 2\alpha - \lambda +1)}{1\cdot 2\cdot 3\cdots \alpha \bullet 1 \cdot 2 \cdot 3 \cdots (\lambda + \alpha)}
\]
ist. Wir wollen daher den folgenden Term
\[
	\ebinom{n}{\alpha +1}\ebinom{n-\alpha -1}{\lambda + \alpha +1} = \Phi
\]
entwickeln, woher
\[
	\Phi = \frac{n(n-1)(n-2) \cdots (n-2\alpha - \lambda -1)}{1\cdot 2\cdot 3 \cdots (\alpha +1) \bullet 1\cdot 2\cdot 3 \cdots (\lambda + \alpha +1)}
\]
wird. Daher berechnen wir also
\[
	\frac{\Phi}{\Pi} = \frac{(n-2\alpha -\lambda)(n -2\alpha - \lambda -1)}{(\alpha +1)(\lambda + \alpha + 1)}
\]
und daher ist
\[
\Pi = \frac{(\alpha +1)(\lambda +\alpha +1)\Phi}{(n-2\alpha -\lambda)(n-2\alpha - \lambda -1)}
\]
\paragraph{§33}
Wir wollen gleich für die folgenden Werte $z'$ und $z''$ die Werte, die $\Phi$ entsprechen, durch $\Phi'$ und $\Phi''$ bezeichnen; weil ja diese aus dem Wert entstehen, wenn $(n+1)$ und $(n+2)$ anstelle von $n$ geschrieben wird, wird nach Ausführung der Entwicklung
\[
	\Phi = \frac{(n+1)n(n-1) \cdots (n-2\alpha - \lambda)}{1\cdot 2\cdot 3 \cdots (\alpha +1) \bullet 1\cdot 2\cdot 3 \cdots (\alpha + \lambda + 1)}
\]
sein, woher klar ist, dass
\[
	\frac{\Phi'}{\Phi} = \frac{n+1}{n-2\alpha - \lambda -1}
\]
sein wird und daher
\[
	\Phi' = \frac{(n+1)\Phi}{n-2\alpha - \lambda -1}.
\]
Auf dieselbe Weise wird
\[
	\Phi'' = \frac{n+2}{n-2\alpha - \lambda}\Phi' = \frac{(n+2)(n+1)\Phi}{(n-2\alpha - \lambda)(n-2\alpha - \lambda - 1)}
\]
sein.
\paragraph{§34}
Daher wollen wir genauso wie oben diesen Ausdruck bilden:
\[
	A\Phi + \frac{B}{n+1}\Phi' + \frac{C}{(n+2)(n+1)}\Phi'',
\]
dessen Wert durch $\Phi$ so ausgedrückt wird:
\[
	\Phi\left( A + \frac{B}{n-2\alpha - \lambda -1)} + \frac{C}{(n-2\alpha -\lambda)(n-2\alpha - \lambda -1)}\right),
\]
wo wiederum die Buchstaben $A,\, B,\, C$ so bestimmt werden müssen, dass die Formel dem vorhergehenden Term $\Pi$ gleich wird. Nachdem also anstelle von $\Pi$ der zuvor durch $\Phi$ ausgedrückte Wert eingesetzt wurde, werden wir diese schon von Brüchen befreite folgende Gleichung erhalten:
\[
A(n-2\alpha - \lambda -1)(n-2\alpha-\lambda) + B(n-2\alpha - \lambda) + C = (\alpha +1)(\alpha + \lambda +1)
\]
\paragraph{§35}
Nachdem also die Entwicklung gemacht wurde und zuerst auf beiden Seiten die Terme gleichgesetzt wurde, die $\alpha\alpha$ involvieren, geht diese Gleichung für die Bestimmung des Buchstaben $A$ hervor:
\[
	4\alpha\alpha A = \alpha\alpha \quad \text{und daher} A = \tfrac{1}{4}.
\]
Wenn auf dieselbe Weise die Terme, die den einfachen Buchstaben $\alpha$ involvieren, gleichgesetzt werden, werden wir zur folgenden Gleichung geführt:
\[
	(4\alpha\lambda - 4n\alpha + 2\alpha)A - 2\alpha B = (\lambda + 2)\alpha,
\]
woher man
\[
	B = -\frac{2n-3}{4}
\]
berechnet. Schließlich geht, nachdem die von $\alpha$ freien Terme gleichgesetzt wurden, diese Gleichung hervor:
\[
\frac{nn-2n\lambda -n + \lambda\lambda + \lambda}{4} - \frac{(n-\lambda)(2n+3)}{4} + C = \lambda +1,
\]
woher
\[
	C = \frac{(n+2)^2}{4} - \frac{\lambda\lambda}{4}
\]
wird.
\paragraph{§36}
Nachdem also die Werte für die einzelnen Terme gefunden wurden, wird immer
\[
	A\Phi + \frac{B}{n+1}\Phi' + \frac{C}{(n+2)(n+1)}\Phi'' = \Pi
\]
sein. Wenn wir daher diese Formel berechnen:
\[
	Az + \frac{B}{n+1}z' + \frac{C}{(n+2)(n+1)}z'',
\]
wird aus den ersten für $\Phi$ angenommenen Terme der vorhergehende der Reihe $z$ entstehen, der $0$ ist; aus den zweiten für $\Phi$ angenommenen Termen wird aber der erste Term $\ebinom{n}{\lambda}$ entstehen; aus den dritten Termen wird der zweite Term $\ebinom{n}{1}\ebinom{n-1}{\lambda +1}$ gebildet; aus den vierten für $\Phi$ angenommenen Term wird der dritte gebildet, welcher $\ebinom{n}{2}\ebinom{n-2}{\lambda +2}$ ist, und so weiter; nachdem diese berechnet worden sind, entsteht diese für $z$ gegebene Reihe selbst
\[
	z = \ebinom{n}{\lambda} + \ebinom{n}{1}\ebinom{n-1}{\lambda +1} + \ebinom{n}{2}\ebinom{n-2}{\lambda +2} + \ebinom{n}{3}\ebinom{n-3}{\lambda +3} + \mathrm{etc.}
\]
Die Relation zwischen $z,\, z',\, z''$ wird also
\[
	Az + \frac{B}{n+1}z' + \frac{C}{(n+2)(n+1)}z'' = z
\]
sein.
\paragraph{§37}
Wir wollen nun anstelle der Buchstaben $A,\,B,\,C$ die gerade gefundenen Werte einsetzen, und die Gleichung zwischen diesen drei Buchstaben wird
\[
	\tfrac{1}{4}z - \frac{2n+3}{4(n+1)}z' + \frac{(n+2)^2 - \lambda\lambda}{4(n+2)(n+1)}z'' = z
\]
sein, die auf diese Form zurückgeführt wird:
\[
	\frac{(n+2)^2 - \lambda\lambda}{(n+2)(n+1)}z'' = \frac{2n+3}{n+1}z' + 3z,
\]
woher man
\[
	z'' = \frac{n+2}{(n+2)^2 - \lambda\lambda}((2n-3)z' + 3(n+1)z)
\]
berechnet.
\paragraph{§38}
Wir wollen nun dem Buchstaben $\lambda$ nacheinander die Werte $0,\, 1,\, 2,\, 3,\, 4,\, \mathrm{etc}$ zuteilen und wir werden die folgenden Relationen für die einzelnen Buchstaben finden:
\begin{align*}
	\frac{(n+2)^2 - 0^2}{(n+2)(n+1)}p'' &= \frac{2n+3}{n+1}p' + 3p \\
	\frac{(n+2)^2 - 1^2}{(n+2)(n+1)}q'' &= \frac{2n+3}{n+1}q' + 3q \\
	\frac{(n+2)^2 - 2^2}{(n+2)(n+1)}r'' &= \frac{2n+3}{n+1}r' + 3r \\
	\frac{(n+2)^2 - 3^2}{(n+2)(n+1)}s'' &= \frac{2n+3}{n+1}s' + 3s \\
	\mathrm{etc.}
\end{align*}
\paragraph{§39}
Weil wir also für den Buchstaben $q$ diese Gleichung haben:
\[
	q'' = \frac{n+2}{(n+1)(n+3)} ((2n+3)q' +3(n+1)q).
\]
Im Fall $n=0$ wird $q=0$ und $q'=1$ sein, woher
\[
	q'' = \tfrac{2}{3}(3\cdot 1 + 3 \cdot 0) = 2
\]
wird. Nun wird für $n=1$ wegen $q=1$ und $q'=2$ 
\[
	q'' = \tfrac{3}{2\cdot 4}(5\cdot 2 + 6 \cdot 1) = 6
\]
sein. Dann wird für den Fall $n=2$ wegen $q=2$ und $q'=6$
\[
	q'' = \tfrac{4}{3\cdot 5}(7\cdot 6 + 9\cdot 2) = 16
\]
sein. Für $n=3$ genommen wird wegen $q=6$ und $q'=16$
\[
	q'' = \tfrac{5}{4\cdot 6} (9 \cdot 16 + 12 \cdot 6) = 45
\]
sein. Aber für den Fall $n=4$ wird wegen $q=16$ und $q'=45$
\[
	q'' = \tfrac{6}{5\cdot 7}(11\cdot 45 + 15\cdot 16) = 126
\]
sein.
\paragraph{§40}
Diese Rechnung ist aber um Vieles arbeitsaufwendiger und ekelhafter als die vorhergehende, die für die Werte des Buchstaben $p$ ausgelegte. Aber es kann eine andere um Vieles leichtere Methode gefunden werden, durch die sich alle Buchstaben $q,\, r,\, s,\, \mathrm{etc}$ durch den Buchstaben $p$ allein mit seinen Derivaten $p',\, p'',\, \mathrm{etc}$ bestimmen lassen werden; dann können daher nämlich, nachdem die Reihe der Zahlen schon hinreichend weit berechnet worden war, die Werte der Buchstaben mit viel geringerem Aufwand berechnet werden; das werden wir im folgenden Abschnitt zeigen.
\section*{\boldmath Bestimmung der Buchstaben $q,\, r,\, s,\, t,\, \mathrm{etc}$ durch den ersten $p$ allein mit seinen Derivaten\unboldmath}
\paragraph{§41}
Nachdem der Kürze wegen unser Trinom
\[
	(1 +x +xx)^n = \chi
\]
gesetzt wurde, wollen wir seine beiden Potenzen $\chi^n$ und $\chi^{n+1}$ so aufteilen, dass die geraden Potenzen von $x$ untereinander geschrieben auftauchen, und zwar auf diese Weise:
\begin{align*}
	\chi^n &= 1 + nx + \cdots + qx^{n-1} + px^n + qx^{n+1} + rx^{n+2} + sx^{n+3} + \mathrm{etc.} \\
	\chi^{n+1} &= 1 + (n+1)x + \cdots + r'x^{n-1} + q'x^n + p'x^{n+1} + q'x^{n+2} + r'x^{n+3} + \mathrm{etc.}
\end{align*}
Nachdem das gemacht wurde, haben wir schon oben bemerkt, dass ein gewisser Koeffizient der unteren Reihe dem unteren mit den beiden vorhergehenden gleich wird.
\paragraph{§42}
Durch dieses Gesetz also werden wir die folgenden Gleichheiten erhalten:
\begin{align*}
	p' &= q + p + q = 2q + p \\
	q' &= r + q + p = p \\
	r' &= s + r + q \\
	\mathrm{etc.}
\end{align*}
woher wir die folgenden Bestimmungen berechnen:
\[
	q = \frac{p'-p}{2},\quad r = q'-q-p,\quad s = r'-r-q,\quad t = s'-s-r,\quad \mathrm{etc.}
\]
\paragraph{§43}
Es ist klar, dass hier die Formel $p'-p$ den Zuwachs der Größe $p$ ausdrückt, während der Exponent $n$ um die Einheit vermehrt wird; weil das durch $\Delta p$ ausgedrückt zu werden pflegt, werden die gefundenen Gleichheiten auf die folgende Weise kürzer beschafft werden können:
\[
	q = \tfrac{1}{2}\Delta p \quad \text{oder}\quad 2q = \Delta p,\quad 2r = 2\Delta q - 2p,\quad 2s = 2\Delta r - 2q,\quad \mathrm{etc.}
\]
\paragraph{§44}
Nachdem aber für diese Differenz der Charakter $\Delta$ benutzt wurde, wird, weil $2q=\Delta p$ ist, $2\Delta q = \Delta\Delta p$ sein und daher
\[
	2r = \Delta\Delta p -2p \quad \text{und daher} \quad 2\Delta r = \Delta^3 p-2\Delta p,
\]
woraus weiter
\[
	2s = \Delta^3 p -3\Delta p \quad \text{und} \quad 2\Delta s = \Delta^4 p -3\Delta\Delta p
\]
wird, also
\[
	2t = \Delta^4 p - 4\Delta\Delta p + 2p \quad \text{und daher} \quad 2\Delta t = \Delta^5 - 4\Delta^3 p + 2\Delta p.
\]
Daher wird weiter
\[
	2u = \Delta^5 p - 5\Delta^3 p + 5\Delta p \quad \text{und daher} \quad 2\Delta u = \Delta^6 p -5\Delta^4 p + 5\Delta\Delta p,
\]
woher man
\[
	2v = \Delta^6 p - 6\Delta^4 p + 9\Delta\Delta p -2p
\]
folgert und so weiter.
\paragraph{§45}
Wenn wir also die numerischen Koeffizienten aufmerksamer betrachten, wird das Bildungsgesetz der Progression entdeckt mit einer den Mathematikern hinreichend bekannten Reihe übereinzustimmen, woher wir für den Wert $z$, für welchen der Index der Ordnung $\lambda$ gesetzt worden ist, die folgende Form erhalten:
\begin{align*}
	2z &= \Delta^{\lambda} p - \lambda\Delta^{\lambda -2} p + \frac{\lambda(\lambda -3)}{1\cdot 2}\Delta^{\lambda -4} p - \frac{\lambda(\lambda -4)(\lambda -5)}{1\cdot 2\cdot 3}\Delta^{\lambda -6}p \\
	&+ \frac{\lambda(\lambda -5)(\lambda -6)(\lambda -7)}{1\cdot 2\cdot 3\cdot 4}\Delta^{\lambda -8}p - \frac{\lambda(\lambda -6)(\lambda -7)(\lambda-8)(\lambda -9)}{1\cdot 2\cdot 3\cdot 4\cdot 5}\Delta^{\lambda -10}p\\
	 &+ \mathrm{etc.}
\end{align*}
welche Reihe nur bis dorthin fortgesetzt werden muss, bis die Indizes von $\Delta$ nicht negativ werden. Wenn wir so $\lambda = 6$ nehmen, in welchem Fall $z = v$ wird, geht aus diesem allgemeinen Bildungsgesetz jedenfalls
\[
	2v = \Delta^6 p - 6\Delta^4 p + 9\Delta^2 p -2p
\]
hervor.
\paragraph{§46}
Damit die Gestalt dieser Reihe besser durchschaut wird, muss man sich erinnern, dass diese Form
\[
	\frac{(x + \sqrt{xx-4})^n}{2^n} + \frac{(x-\sqrt{xx-4})^n}{2^n}
\]
in die folgende Reihe aufgelöst wird:
\[
	x^n - nx^{n-2} + \frac{n(n-3)}{1\cdot 2}x^{n-4} - \frac{n(n-4)(n-5)}{1\cdot 2\cdot 3}x^{n-6} + \mathrm{etc.}
\]
Auf diese Weise also ist unserem Anliegen schon genügt worden, weil wir alle Buchstaben $q,~ r,~ s,~ t,~ \mathrm{etc}$ durch den ersten $p$ allein und seine Derivaten $p',~ p'',~ p''',~ \mathrm{etc}$ ausgedrückt gefunden haben.
\section*{\boldmath Bestimmung der Größe $p$ durch eine endliche Integralformel \unboldmath}
\paragraph{§47}
Weil durch die dritte oben [§16] erörterte Formel
\[
	p = 1 + \ebinom{2}{1}\ebinom{n}{2} +\ebinom{4}{2}\ebinom{n}{4} + \ebinom{6}{3}\ebinom{n}{6} + \mathrm{etc.}
\]
ist, wird der beliebige Term im Allgemeinen $\ebinom{2\alpha}{\alpha}\ebinom{n}{2\alpha}$ sein, dem dieser $\ebinom{2\alpha +2}{\alpha +1}\ebinom{n}{2\alpha +2}$ folgt.

Weil also, nachdem die Entwicklung gemacht worden ist,
\[
	\ebinom{2\alpha}{\alpha} = \frac{2\alpha(2\alpha -1)(2\alpha -2) \cdots (\alpha +1)}{1\cdot 2\cdot 3 \cdots \alpha}
\]
ist, und auf ähnliche Weise wird
\[
	\ebinom{2\alpha +2}{\alpha +1} = \frac{(2\alpha +2)(2\alpha +1)2\alpha \cdots (\alpha + 2)}{1\cdot 2 \cdot 3 \cdots (\alpha + 1)}
\]
sein; diese letzte Form gibt durch die erste geteilt den Quotienten 
\[
	\frac{(2\alpha + 2)(2\alpha +1)}{(\alpha +1)^2} = \frac{2(2\alpha +1)}{\alpha +1}
\]
und so wird
\[
	\ebinom{2\alpha +2}{\alpha +1} = \frac{4\alpha + 2}{\alpha +1}\ebinom{2\alpha}{\alpha}
\]
sein.
\paragraph{§48}
Nachdem also diese Reduktion angewendet worden ist, wird für $\alpha = 1 $ 
\[
\ebinom{4}{2} = \frac{6}{2}\ebinom{2}{1}
\]
sein; für $\alpha = 2$ wird
\[
\ebinom{6}{3} = \frac{10}{3}\ebinom{4}{2} = \frac{10}{3}\cdot\frac{6}{2}\cdot\frac{2}{1}
\]
sein; wenn $\alpha = 3$ ist, wird
\[
\ebinom{8}{4} = \frac{14}{4}\ebinom{6}{3} = \frac{14}{4}\cdot\frac{10}{3}\cdot\frac{6}{2}\cdot\frac{2}{1}
\]
sein; wenn dann	$\alpha = 4$ ist, wird
\[
\ebinom{10}{5} = \frac{18}{5}\ebinom{8}{4} = \frac{18}{5}\cdot\frac{14}{4}\cdot\frac{10}{3}\cdot\frac{6}{2}\cdot\frac{2}{1}
\]
sein, und so weiter. Nachdem also diese Werte eingeführt wurden, wird durch gewöhnliche numerische Faktoren
\begin{align*}
p &= 1 + \frac{2}{1}\ebinom{n}{2} + \frac{2\cdot 6}{1\cdot 2}\ebinom{n}{4} + \frac{2\cdot 6\cdot 10}{1\cdot 2\cdot 3}\ebinom{n}{6} + \frac{2\cdot 6\cdot 10\cdot 14}{1\cdot 2\cdot 3\cdot 4}\ebinom{n}{8} \\
&+ \frac{2\cdot 6\cdot 10\cdot 14 \cdot 18}{1\cdot 2\cdot 3\cdot 4\cdot 5}\ebinom{n}{10} + \mathrm{etc.}
\end{align*}
sein.
\paragraph{§49}
Nun wollen wir also sehen, wie diese endliche Integralformel gefunden werden muss, deren Integral zwischen gegebenen Grenzen eingeschlossen zu dieser Reihe selbst führt. Für dieses Ziel muss diese Formel betrachtet werden, $(1+x)^n$, deren Entwicklung ja diese Reihe liefert:
\[
	1 + \ebinom{n}{1}x + \ebinom{n}{2}xx + \ebinom{n}{3}x^3 + \mathrm{etc,}
\]
deren Terme wechselweise schon unsere Charaktere des Buchstaben n enthalten.
\paragraph{§50}
Wir wollen daher diese Reihe in zwei Teile teilen, nach den abwechselnden Termen, und wollen
\[
	M = 1 + \ebinom{n}{2}xx + \ebinom{n}{4}x^4 + \ebinom{n}{6}x^6 + \mathrm{etc.}
\]
\[
	N = \ebinom{n}{1}x + \ebinom{n}{3}x^3 + \ebinom{n}{5}x^5 + \ebinom{n}{7}x^7 + \mathrm{etc.}
\]
setzen, so dass
\[
	(1+x)^n = M + N
\]
ist. Nun wollen wir aber untersuchen, wie die erste Reihe, $M$, durch analytische Operationen behandelt werden muss, dass die vorgelegte Reihe selbst oder der Wert von $p$ daraus entspringt.
\paragraph{§51}
Um das zu bewirken, wollen wir die Größe $M$ mit einem gewissen Differential $\partial \nu$ einer Funktion von $x$ multiplizieren, und die folgenden Integrationen so bestimmen, dass sie zwischen bestimmten Grenzen, wie z.\,B. $x=a$ bis hin zu $x=b$, eingeschlossen werden, welche Bedingungen so beschaffen sein müssen, dass es den folgenden Bedingungen genügt:
\begin{align*}
	1.) \int xx\partial \nu &= \frac{2}{1}\nu \\
	2.) \int x^4\partial \nu &= \frac{2\cdot 6}{1\cdot 2}\nu \\
	3.) \int x^6\partial \nu &= \frac{2\cdot 6\cdot 10}{1\cdot 2\cdot 3}\nu \\
	4.) \int x^8\partial \nu &= \frac{2\cdot 6\cdot 10\cdot 14}{1\cdot 2\cdot 3\cdot 4}\nu \\
	\mathrm{etc.}
\end{align*}
Auf diese Weise nämlich wird das Integral $\int M\partial \nu$ diese Reihe ergeben:
\[
	\nu + \frac{2}{1}\ebinom{n}{2}{\nu} + \frac{2\cdot 6}{1\cdot 2}\ebinom{n}{4}{\nu} + \frac{2\cdot 6\cdot 10}{1\cdot 2\cdot 3}\ebinom{n}{6}{\nu} + \mathrm{etc,}
\]
sodass wir auf diese Weise das, was wir sehen, erreichen,
\[
	p = \int\frac{M\partial \nu}{\nu}
\]
\paragraph{§52}
Eine beliebige dieser Integralformeln, die wir hier erörtert haben, hängt so von der vorhergehenden ab, dass
\begin{align*}
	\int xx \partial \nu &= \frac{2}{1}\int\partial\nu \\
	\int x^4\partial \nu &= \frac{6}{2}\int xx \partial \nu \\
	\int x^6 \partial \nu &= \frac{10}{3}\int x^4 \partial \nu \\
	\int x^8 \partial \nu &= \frac{14}{4}\int x^6 \partial \nu \\
	\mathrm{etc.}
\end{align*}
ist und so muss im Allgemeinen bewirkt werden, dass
\[
	\int x^{2m}\partial \nu = \frac{4m-2}{m}\int x^{2m-2}\partial\nu
\]
wird. Diese Reduktion muss natürlich Geltung haben, nachdem die Integrale zwischen vorgeschriebenen Grenzen angenommen worden sind, natürlich von $x=a$ bis hin zu $x=b$, welche Terme zwar noch nicht bekannt sind, aber an die vorgeschriebene Bedingung selbst angepasst werden müssen.
\paragraph{§53}
Weil also für diese Integrationsgrenzen
\[
	m\int x^{2m} \partial \nu = (4m-2)\int x^{2m-2}\partial \nu
\]
sein muss, wollen wir setzen, dass im Allgemeinen
\[
	m\int x^{2m}\partial \nu = (4m-2)\int x^{2m-2}\partial \nu + \Pi x^{2m-1}
\]
ist, wo natürlich $\Pi$ eine Funktion solcher Art ist, dass der angeknüpfte Teil $\Pi x^{2m-1}$ für jede der beiden Grenzen, wie für $x=a$ so für $x=b$, verschwindet. Diese Gleichung gibt differentiert und durch $x^{2m-2}$ geteilt
\[
	mxx\partial\nu = (4m-2)\partial\nu + (2m-1)\Pi \partial x + x\partial\Pi,
\]
welche Gleichung für alle numerischen Werte $m$ gelten muss.
\paragraph{§54}
Daher also wird diese Gleichung in zwei Teile geteilt werden müssen, von denen die eine einzig mit dem Buchstaben $m$ versehene Terme enthalte, die andere aber die übrigen, welche zwei Gleichungen, also
\[
	xx\partial\nu = 4\partial \nu + 2\Pi\partial x
\]
und
\[
	0 = -2\partial \nu - \Pi\partial x + x \partial \Pi,
\]
sein werden. Aus der ersten wird
\[
	\partial \nu = \frac{2 \Pi \partial x}{xx-4};
\]
aus der anderen aber wird
\[
	\partial \nu = \frac{x\partial\Pi - \Pi\partial x}{2},
\]
welche beiden Werte einander gleichgesetzt diese Gleichung liefern:
\[
	4\Pi\partial x = (xx-4)(x\partial\Pi - \Pi\partial x) = x^3\partial\Pi - xx\Pi\partial x - 4x\partial \Pi + 4\Pi\partial x
\]
und daher berechnet man
\[
	\frac{\partial \Pi}{\Pi} = \frac{x\partial x}{xx-4},
\]
woher durch Integrieren
\[
	\log{\Pi} = \log{\sqrt{xx-4}}
\]
wird und daher
\[
	\Pi = C\sqrt{xx-4}
\]
oder auch
\[
	\Pi = C\sqrt{4-xx};
\]
nachdem dieser Wert gefunden wurde, erhalten wir unser angenommenes Differential
\[	
	\partial \nu = \frac{2C\partial x}{\sqrt{4-xx}},
\]
woher
\[
	\nu = 2C\arcsin{\left(\frac{x}{2}\right)}
\]
wird.
\paragraph{§55}
Wir wollen nun die angefügte Formel
\[
	\Pi x^{2m-1} = Cx^{2m-1}\sqrt{4-xx}
\]
betrachten, welche wir entdecken auf dreifache Weise verschwinden zu können: zuerst natürlich, wann immer $x=0$ ist, ausgenommen im Fall $m=0$; zweitens im Fall, in dem $x=2$ ist; und drittens in dem Fall, in dem $x=-2$ ist, aus welchen also beide Integrationsgrenzen $a$ und $b$ ausgesucht werden müssen. Diese beiden Integrationsgrenzen müssen aber so gewählt werden, dass auch der andere Teil der Integration, $\int N\partial x$, angenehm ausgedrückt wird. Weil wir nämlich
\[
	(1+x)^n = M + N
\] 
gesetzt haben, ist auch auf das Integral $\int N \partial \nu$ zu achten; wenn dieses völlig verschwinden würde, wäre das für die Integrationsgrenzen ohne Zweifel sehr angenehm, dann wäre nämlich 
\[
	\int (M+N) \partial \nu
\]
oder
\[
	\int \partial \nu (1+x)^n = \int M \partial\nu;
\]
als logische Konsequenz würden wir
\[
	p = \frac{\int M \partial \nu}{\nu}
\]
haben.
\paragraph{§56}
Oben haben wir aber
\[
	N = \ebinom{n}{1}x + \ebinom{n}{3}x^3 + \ebinom{n}{5}x^5 + \ebinom{n}{7}x^7 + \mathrm{etc.}
\]
gesetzt, woher
\[
	\int N \partial \nu = \ebinom{n}{1}\int x \partial \nu + \ebinom{n}{3} \int x^3 \partial \nu + \ebinom{n}{5}\int x^5 \partial \nu + \mathrm{etc.}
\]
gebildet wird, wodurch dieselben Reduktionen, welche wir für den Buchstaben $M$ angestellt haben, eine beliebige Integralformel zu der vorhergehenden mithilfe der Reduktion
\[
	\int x^{2m}\partial \nu = \frac{4m-2}{m}\int x^{2m-2}\partial \nu
\]
zurückgeführt werden kann. Für $m=\frac{3}{2}$ genommen wird nämlich
\[
	\int x^3 \partial \nu = \frac{8}{3}\int \partial \nu
\]
sein. Für $m=\frac{5}{2}$ genommen wird
\[
	\int x^5 \partial \nu = \frac{16}{5}\int x^3 \partial \nu
\]
sein. Für $m=\frac{7}{2}$ genommen wird
\[
	\int x^7 \partial \nu = \frac{24}{7}\int x^5 \partial \nu
\]
sein $\mathrm{etc}$, woher klar ist, wenn nur $\int x \partial \nu$ verschwände, dass dann auch alle folgenden verschwinden werden.
\paragraph{§57}
Weil wir also
\[
	\partial \nu = \frac{2C\partial x}{\sqrt{4-xx}}
\]
gefunden haben, wird
\[
	x\partial \nu = \frac{2Cx\partial x}{\sqrt{4-xx}}
\]
sein und daher
\[
	\int x \partial \nu = 2C\sqrt{4-xx},
\]
welcher Ausdruck in beiden Fällen entweder bei $x=+2$ oder $x=-2$ verschwindet. Wenn wir deshalb die Integrationsgrenzen $x=2$ und $x=-2$ setzen, verschwinden nicht nur jene angefügte Teile $\Pi x^{2m-1}$, sondern auch der ganze Wert des Integrals $\int N \partial \nu$, und daher genügen wir in diesem Fall unserer Frage vollkommen, weil
\[
	p = \frac{\int \partial \nu (1+x)^n}{\nu}
\]
ist.
\paragraph{§58}
Weil wir also
\[
	\partial \nu = \frac{2C\partial x}{\sqrt{4-xx}}
\]
gefunden haben, dessen Integral so genommen wird, dass es für $x=2$ gesetzt verschwindet, wird
\[
	\nu = 2C\arcsin{\frac{x}{2}} - 2C\frac{\pi}{2}
\]
sein, welcher Ausdruck auf diesen zurückgeführt wird:
\[
	\nu = -2C\arccos{\frac{x}{2}};
\]
daher geht, nachdem das Integral bis hin zur anderen Grenze $x=-2$ erstreckt wurde, $\nu = -2C\pi$ hervor. Nach Einsetzen dieser Werte wird also die gesuchte Formel
\[
	p = -\frac{1}{\pi}\int\frac{(1+x)^n \partial x}{\sqrt{4-xx}}
\]
sein. Diese Integralformel wird natürlich von der Grenze $x=2$ bis hin zur Grenze $x=-2$ erstreckt den wahren Wert von $p$ liefern.
\paragraph{§59}
Damit wir diese Formel vereinfachen, wollen wir $x=2\cos{\varphi}$ setzen, wo klar ist, dass im Fall $x=2$ der Winkel $\varphi = 0$ wird; im Fall $x=-2$ aber $\varphi = \pi$, sodass nach Einführung dieses Winkels das Integral von der Grenze $\varphi = 0$ bis hin zu $\varphi = \pi$ genommen werden muss; dann aber wird
\[
	\partial x = -2\partial \varphi \sin{\varphi} \quad \text{und} \quad \sqrt{4-xx} = 2\sin{\varphi}
\]
sein, nach welcher Substitution wir diese Gleichung erhalten werden:
\[
	p = +\frac{1}{\pi}\int (1+2\cos{\varphi})^n \partial \varphi \qquad \begin{bmatrix}
	\text{von}~ \varphi = 0 \\
	\text{bis}~ \varphi = \pi
	\end{bmatrix}
\]
\section*{Bestimmung der übrigen Buchstaben durch endliche Integralformeln}
\paragraph{§60}
Das kann leicht durch die Relationen, die wir oben zwischen diesen Buchstaben angegeben haben, geleistet werden. Zuerst hatten wir natürlich $2q = \Delta p = p'-p$, wo $p'$ aus $p$ entsteht, wenn $n+1$ anstelle von $n$ geschrieben wird. Weil wir ja gerade
\[
	p = \frac{1}{\pi}\int (1 + 2\cos{\varphi})^n \partial \varphi
\]
gefunden haben, wird
\[
	p' = \frac{1}{\pi}\int (1 + 2\cos{\varphi})^{n+1} \partial \varphi
\]
sein und daher wird also
\[
	p' - p = \frac{2}{\pi}\int \cos{\varphi} (1 +2\cos{\varphi})^n \partial \varphi
\]
sein, nach Einsetzen welchen Wertes man
\[
	q = \frac{1}{\pi}\int\partial \varphi \cos{\varphi} (1 + 2\cos{\varphi})^n \qquad \begin{bmatrix}
	\text{von}~ \varphi = 0 \\
	\text{bis}~ \varphi = \pi
	\end{bmatrix}
\]
finden wird; daher wird also
\[
	q' = \frac{1}{\pi} \int \partial \varphi \cos{\varphi} (1 + 2\cos{\varphi})^{n+1}
\]
sein.
\paragraph{§61}
Oben aber haben wir gesehen, dass $r = q' - q - p$ ist, nun wird aber
\[
	q'-q = \frac{2}{\pi}\int \partial \varphi \cos^2{\varphi} (1 + 2\cos{\varphi})^n
\]
sein. Wenn daher also $p$ abgezogen wird, finden wir wegen $2\cos^2{\varphi} - 1 = \cos{2\varphi}$ den Buchstaben
\[
	r = \frac{1}{\pi}\int \partial \varphi \cos{2\varphi} (1 + 2\cos{\varphi})^n,
\]
woher wiederum
\[
	r' = \frac{1}{\pi}\int \partial \varphi \cos{2\varphi}(1 + 2\cos{\varphi})^{n+1}
\]
wird.
\paragraph{§62}
Weil wir ja oben $s = r' - r - q$ gefunden haben, werden wir hier zuerst
\[
	r' - r = \frac{2}{\pi}\int \partial\varphi \cos{\varphi}\cos{2\varphi}(1 + 2\cos{\varphi})^n
\]
haben. Wenn daher also $q$ abgezogen wird, wird wegen $2\cos{\varphi}\cos{2\varphi} - \cos{\varphi} = \cos{3\varphi}$
\[
	s = \frac{1}{\pi}\int\partial\varphi\cos{3\varphi}(1+2\cos{\varphi})^n
\]
sein. Auf ähnliche Weise ist schon klar, dass
\[
	t = \frac{1}{\pi}\int\partial\varphi\cos{4\varphi} (1 + 2\cos{\varphi})^n
\]
sein wird, auf dieselbe Weise wird man finden, dass
\[
	u = \frac{1}{\pi}\int\partial\varphi\cos{5\varphi}(1+2\cos{\varphi})^n
\]
sein wird und daher wird im Allgemeinen
\[
	z = \frac{1}{\pi}\int\partial\varphi\cos{\lambda\varphi}(1 + 2\cos{\varphi})^n
\]
sein.
\paragraph{§63}
Weil ja die Analysis, die wir benutzt haben, völlig einzigartig und ungewöhnlich ist, wird es nicht unpassend sein, dass die Gültigkeit dieser Formeln durch einen analytischen Beweis untermauert wird, welcher sich für die einzelnen quasi mit einem Schritt auf die folgende Weise führen lässt. Es wird von der Entwicklung der Formel $(1+2\cos{\varphi})^n$ aus zu beginnen sein, die auf diese Reihe führt:
\[
	1 + \ebinom{n}{1}2\cos{\varphi} + \ebinom{n}{2}4\cos^2{\varphi} + \ebinom{n}{3}8\cos^3{\varphi} + \ebinom{n}{4}16\cos^4{\varphi} + \mathrm{etc.}
\]
Durch bekannte Reduktionen der Winkel ist aber bekannt, dass
\begin{align*}
2\cos{\varphi} &= 2\cos{\varphi} \\
4\cos^2{\varphi} &= 2\cos{2\varphi} + 2 \\
8\cos^3{\varphi} &= 2\cos{3\varphi} + 6\cos{\varphi} \\
16\cos^4{\varphi} &= 2\cos{4\varphi} + 8\cos{2\varphi} + 6\\
32\cos^5{\varphi} &= 2\cos{5\varphi} + 10\cos{3\varphi} + 20\cos{\varphi} \\
\vdots \\
2^{\alpha}\cos^{\alpha}{\varphi} &= 2\cos{\alpha\varphi} + 2\ebinom{\alpha}{1}\cos{(\alpha - 2)}\varphi + 2\ebinom{\alpha}{2}\cos{(\alpha - 4)}\varphi \\
&+ 2\ebinom{\alpha}{3}\cos{(\alpha - 6)}\varphi + \mathrm{etc.}
\end{align*}
sein wird; es ist natürlich zu bemerken, dass, wann immer der letzte Term konstant ist, er nur einfach genommen werden muss; außerdem müssen aber auch die Cosinus der negativen Winkel völlig weggelassen werden.
\paragraph{§64}
Nachdem diese also auf gewohnte Weise aufgeteilt worden sind, wird
\begin{align*}
	(1+2\cos{\varphi})^n &= 1 + \ebinom{n}{1}2\cos{\varphi} + \ebinom{n}{2}2(\cos{2\varphi + 1}) + 2\ebinom{n}{3}(\cos{3\varphi} + 3\cos{\varphi}) \\
	&+ 2\ebinom{n}{4}(\cos{4\varphi} + 4\cos{2\varphi} + 3) + 2\ebinom{n}{5}(\cos{5\varphi} + 5\cos{3\varphi} + 10\cos{\varphi}) \\
	&+ 2\ebinom{n}{6}(\cos{6\varphi} + 6\cos{4\varphi} + 15\cos{2\varphi} + 10) + \mathrm{etc.}
\end{align*}
sein, woher die folgenden Integrationen zu suchen sind.
\paragraph{§65}
Wir wollen vom ersten Buchstaben $p$ aus beginnen, wo diese Reihe mit $\partial \varphi$ multipliziert werden und integriert werden muss. Weil daher im Allgemeinen
\[
	\int \partial \varphi \cos{m\varphi} = \frac{1}{m}\sin{m\varphi}
\]
ist, verschwindet dieser Wert für $\varphi = 0$ gesetzt; für die andere Integrationsgrenze $\varphi = \pi$ verschwindet er natürlich, wenn alle Zahlen $n$ ganze sind. Zur Integration also bleiben allein die absoluten Terme übrig; dann aber wird, nachdem das Integral auf gewohnte Weise genommen wurde, $\int\partial \varphi = \pi$ sein, nach Bemerkung wovon unser Integral
\[
	\int\partial \varphi (1 + 2\cos{\varphi})^n = \pi + 2\ebinom{n}{2}\pi + 6\ebinom{n}{4}\pi + 20\ebinom{n}{6}\pi + \mathrm{etc.}
\]
sein wird. Wenn daher hier die allgemeine oben gegebene Form bemüht wird, führe man diese Koeffizienten auf diese Formen $\ebinom{2}{1},\,\ebinom{4}{2},\,\ebinom{6}{3},\,\mathrm{etc}$ zurück, genauso wie es die Gültigkeit der Formel erfordert. Es wird nämlich
\[
	p = \frac{1}{\pi}\int\partial\varphi (1 + 2\cos{\varphi})^n = 1 + \ebinom{2}{1}\ebinom{n}{2} + \ebinom{4}{2}\ebinom{n}{4} + \ebinom{6}{3}\ebinom{n}{6} + \mathrm{etc.}
\]
sein.
\paragraph{§66}
Wir wollen zum zweiten Buchstaben, $q$, vordringen, wo die obere Reihe mit $\partial\varphi\cos{\varphi}$ multipliziert und integriert werden muss. Dabei bemerke man, dass im Allgemeinen
\[
	\int\partial\varphi\cos{\varphi}\cos{m\varphi} = \frac{1}{2(m+1)}\sin{(m+1)\varphi} + \frac{1}{2(m-1)}\sin{(m-1)\varphi}
\]
ist, welcher Ausdruck für $\varphi = \pi$ gesetzt verschwindet, mit Ausnahme des Falles $m=1$, in welchem natürlich
\[
	\int\partial\varphi\cos{\varphi}\cos{\varphi} = \frac{1}{2}\varphi = \frac{\pi}{2}
\]
wird. Daraus sieht man ein, dass aus der oberen Reihe die anderen Terme hier nicht in die Rechnung eingehen, wenn diese nicht $\cos{\varphi}$ enthalten, welche
\[
	2\ebinom{n}{1}\cos{\varphi} + 2\ebinom{3}{1}\ebinom{n}{3}\cos{\varphi} + 2\ebinom{5}{2}\ebinom{n}{5}\cos{\varphi} + 2\ebinom{7}{3}\ebinom{n}{7}\cos{\varphi} + \mathrm{etc.}
\]
sind. Diese Terme werden aber mit $\partial\varphi\cos{\varphi}$ multipliziert und integriert, wegen
\[
	\int 2\partial\varphi\cos^2{\varphi} = \pi
\]
durch $\pi$ geteilt den Wert
\[
	q = \ebinom{n}{1} + \ebinom{3}{1}\ebinom{n}{3} + \ebinom{5}{2}\ebinom{n}{5} + \ebinom{7}{3}\ebinom{n}{7} + \mathrm{etc.}
\]
geben.
\paragraph{§67}
Für den Buchstaben $r$ muss die obere Reihe mit $\partial\varphi\cos{2\varphi}$ multipliziert werden. Weil daher im Allgemeinen
\[
	\cos{2\varphi}\cos{m\varphi} = \tfrac{1}{2}\cos{(m+2)\varphi} + \tfrac{1}{2}\cos{(m-2)\varphi}
\]
ist, verschwindet, indem man mit $\partial\varphi$ multipliziert, das Integral für die Grenze $\varphi = \pi$ immer, einzig ausgenommen im Fall $m=2$, indem natürlich
\[
	\int\partial\varphi\cos^2{2\varphi} = \frac{\pi}{2}
\]
wird. Hier also gehen aus der oberen Reihe einzig die mit $\cos{2\varphi}$ multiplizierten Terme in die Rechnung ein, welche
\[
	2\cos{2\varphi}\left( \ebinom{n}{2} + \ebinom{4}{1}\ebinom{n}{4} + \ebinom{6}{2}\ebinom{n}{6} + \ebinom{8}{3}\ebinom{n}{8} + \mathrm{etc.}\right)
\]
sind. Weil also $2\int\partial\varphi\cos^2{\varphi} = \pi$ ist, findet man, indem man alle Terme sammelt und durch $\pi$ teilt,
\[
	r = \ebinom{n}{2} + \ebinom{4}{1}\ebinom{n}{4} + \ebinom{6}{2}\ebinom{n}{6} + \ebinom{8}{3}\ebinom{n}{8} + \mathrm{etc.}
\]
\paragraph{§68}
Damit das deutlicher herauskommt und leichter auf den allgemeinen Wert $z$ angewendet werden kann, wollen wir die Entwicklung der Potenz $(1 + 2\cos{\varphi})^n$ sofort nach den Kosinus der Vielfachen des Winkels $\varphi$ auf diese Weise aufteilen:
\begin{align*}
&(1+2\cos{\varphi})^n = 1 + \ebinom{2}{1}\ebinom{n}{2} + \ebinom{4}{2}\ebinom{n}{4} + \ebinom{6}{3}\ebinom{n}{6} + \ebinom{8}{4}\ebinom{n}{8} + \mathrm{etc.} \\
&+ 2\cos{\varphi}\left(\ebinom{n}{1} +\ebinom{3}{1}\ebinom{n}{3} +\ebinom{5}{2}\ebinom{n}{5} +\ebinom{7}{3}\ebinom{n}{7} +\ebinom{9}{4}\ebinom{n}{9} +\mathrm{etc.}\right) \\
&+ 2\cos{2\varphi}\left(\ebinom{n}{2} + \ebinom{4}{1}\ebinom{n}{4} + \ebinom{6}{2}\ebinom{n}{6} + \ebinom{8}{3}\ebinom{n}{8} + \ebinom{10}{4}\ebinom{n}{10} + \mathrm{etc.}\right) \\
&+ 2\cos{3\varphi}\left(\ebinom{n}{3} + \ebinom{5}{1}\ebinom{n}{5} + \ebinom{7}{2}\ebinom{n}{7} + \ebinom{9}{3}\ebinom{n}{9} + \ebinom{11}{4}\ebinom{n}{11} + \mathrm{etc.}\right) \\
&\qquad\quad\vdots \\
&+2\cos{\lambda\varphi}\left(\ebinom{n}{\lambda} + \ebinom{\lambda +2}{1}\ebinom{n}{\lambda +2} + \ebinom{\lambda +4}{2}\ebinom{n}{\lambda +4} + \ebinom{\lambda +6}{3}\ebinom{n}{\lambda + 6} + \mathrm{etc.}\right) \\
&+ \mathrm{etc.}
\end{align*}
\paragraph{§69}
Wenn wir also diese Gleichung mit $\partial\varphi\cos{\lambda\varphi}$ multiplizieren und integrieren, werden alle Integrale, die durch vorgeschriebenen Grenzen eingeschlossen wurden, mit Ausnahme des Gliedes $2\cos{\lambda\varphi(\dots)}$, verschwinden, deshalb weil das Produkt $2\cos^2{\lambda\varphi}$ ein konstantes Glied enthält, woher durch Integration $\pi$ entsteht, sodass
\begin{align*}
	\int\partial\cos{\varphi\lambda\varphi}(1+2\cos{\varphi})^n &= \pi\left(\ebinom{n}{\lambda} + \ebinom{\lambda + 2}{1}\ebinom{n}{\lambda + 2}\right. \\
	&+ \left. \ebinom{\lambda + 4}{2}\ebinom{n}{\lambda + 4} + \mathrm{etc.}\right)
\end{align*}
ist, welcher Wert durch $\pi$ geteilt den oben gefundenen Wert von $z$ liefert; daher ist die Gültigkeit dieser neuen Ausdrücke glänzend gezeigt worden.
\paragraph{§70}
Im Übrigen, wenn wir die einzelnen Reihen des vorletzten Termes einmal auf leichte Weise betrachten, entdecken wir, dass sie unseren Buchstaben $p,\, q,\, r,\, s,\, \mathrm{etc}$ selbst gleich sind, sodass nun
\[
	(1+2\cos{\varphi})^n = p + 2q\cos{\varphi} + 2r\cos{2\varphi} + 2s\cos{3\varphi} + 2t\cos{4\varphi} + \mathrm{etc.}
\]
ist, wo nun zugleich der Grund klar ist, warum die Buchstaben $q,\, r,\, s,\, \mathrm{etc}$ verdoppelt werden, was ja darin liegt, dass in der Enticklung der Formel $(1+x+xx)^n$ der Buchstabe $p$ nur einmal in der Mitte, die übrigen Buchstaben aber zweimal, von der Mitte gleich weit entfernt, auftauchen. Daraus ist diese außergewöhnliche Beziehung zwischen jenen beiden Potenzen $(1+x+xx)^n$ und $(1+2\cos{\varphi})^n$ der nächsten Aufmerksamkeit würdig anzusehen.
\section*{\boldmath Untersuchung der Summe der Reihe \\ $P = 1 +x +3xx+ 7x^3 + 19x^4 + \cdots + px^n + p'x^{n+1 } + \mathrm{etc.}$ \unboldmath}
\paragraph{§71}
Weil ja der allgemeine Term dieser Reihe $px^n$ ist, dem $p'x^{n+1}$ und $p''x^{n+2}$ folgen, haben wir zwischen diesen drei Größen $p,\, p',\, p''$ oben [§38] diese Relation gefunden:
\[
	(n+2)p'' = (2n+3)p' + 3(n+1)p,
\]
welche wir zu unserem Nutzen angewandt auf diese Weise darstellen wollen:
\[
	3(n+1)p + (n+1)p' + (n+2)p' - (n+2)p'' = 0
\]
\paragraph{§72}
Weil ja unsere Reihe gleich
\[
	1 + x + 3xx + 7x^3 + 19x^4 + \cdots + px^n + p'x^{n+1} + p''x^{n+2} + \mathrm{etc.}
\]
ist, wollen wir Operationen solcher Art ausführen, mit denen die gerade erwähnte Relation erhalten wird; das wird auf die folgende Weise am besten geschehen:
\begin{align*}
\frac{3\partial Px}{\partial x} &= 3 + 6x + 27xx + \cdots + 3(n+1)px^n + \mathrm{etc.} \\
+\frac{\partial P}{\partial x} &= 1 + 6x + 21xx + \cdots + (n+1)p'x^n + \mathrm{etc.} \\
+\frac{\partial P x}{x\partial x} &= \frac{1}{x} + 2 + 9x + 28xx + \cdots + (n+2)p'x^n + \mathrm{etc.} \\
-\frac{\partial P}{x\partial x} &= -\frac{1}{x} - 6 - 21x - 76xx - \cdots -(n+2)p''x^n + \mathrm{etc.}
\end{align*}
Man fasse diese vier Summen gleich in eine Summe zusammen, und wir werden die folgende Gleichung erhalten:
\[
	\frac{3\partial Px}{\partial x} + \frac{\partial P}{\partial x} + \frac{\partial Px}{x\partial x} - \frac{\partial P }{x\partial x} = 0,
\]
weil sich ja alle Terme gegenseitig aufheben.
\paragraph{§73}
Auf diese Weise also sind wir zu einer endlichen Differentialgleichung ersten Grades geführt worden, welche mit $x\partial x$ multipliziert und geordnet sich so verhalten wird:
\[
	P\partial x (3x + 1) + \partial P(3xx+2x-1) = 0,
\]
woher also
\[
	\frac{\partial P}{P} = \frac{\partial x (1+3x)}{1-2x-3xx}
\]
wird, welche Gleichung integriert
\[
	\log{P} = -\tfrac{1}{2}\log{(1-2x-3xx)} + \log{C}
\]
liefert; als logische Konsequenz ist
\[
	P = \frac{C}{\sqrt{1-2x-3xx}};
\]
um hier die Konstante C zu bestimmen, bemerke man, dass unsere vorgelegte Reihe nur im Fall $x=0$ gleich $P=1$ liefert, woher klar ist, dass $C=1$ genommen werden muss, sodass die Summe der Reihe
\[
	P = \frac{1}{\sqrt{1-2x-3xx}}
\]
ist. 
\paragraph{§74}
Wider der Erwartung also sind wir zu einer algeraischen Summe gelangt, welcher Ausdruck sogar so beschaffen ist, dass er in eine Reihe verwandelt unsere Reihe selbst erzeugt; es wird der Mühe wert sein, das gezeigt zu haben. Weil also
\[
	P = (1-2x-3xx)^{-\frac{1}{2}}
\]
ist, welche es genügen wird, nur bis zur dritten Potenz entwickelt zu haben. Auf diese Weise werden wir
\begin{alignat*}{6}
	P &= 1 &+ x &+ &\frac{3}{2}xx &+ \frac{9}{2}x^3 &+ \mathrm{etc.} \\
	 &     &   &   &+\frac{3}{2}xx &+ \frac{5}{2}x^3 &+ \mathrm{etc.} \\
	  &= 1 &+ x &+ &3xx &+ 7x^3 &+ \mathrm{etc.}
\end{alignat*}
erhalten, welche also vollkommen übereinstimmt.
\paragraph{§75}
Aber in der Tat kann dieselbe Summe noch auf eine andere Weise untersucht werden, natürlich von der Integralformel aus, welche wir für den Wert des Buchstaben $p$ gefunden haben
\[
	p = \frac{1}{\pi}\int \partial\varphi (1+2\cos{\varphi})^n \qquad \begin{bmatrix}
	\text{von}~ \varphi = 0 \\
	\text{bis}~ \varphi = \pi \\
	\end{bmatrix}
\]
Diese Formel nämlich, für $n=0$ genommen, gibt mit $x$ multipliziert den zweiten Term, $x$; im Fall $n=2$ weiter, gibt sie mit $xx$ multipliziert den dritten Term, $3xx$; nachdem das beobachtet wurde, wird die gesuchte Summe so dargestellt werden können:
\[
	P = \frac{1}{\pi}\int\partial\varphi (1 + x(1+\cos{\varphi}) + xx(1+2\cos{\varphi})^2 + x^3(1+2\cos{\varphi})^3 + \mathrm{etc.}),
\]
wo natürlich zu bemerken ist, dass in dieser Integration die Größe $x$ als Konstante betrachtet wird, weil ja einzig der Winkel $\varphi$ variabel ist.
\paragraph{§76}
Es ist aber klar, dass die unendliche Reihe, mit welcher das Element $\partial \varphi$ multipliziert werden muss, eine geometrische ist, deren Summe also
\[
	\frac{1}{1-x(1+2\cos{\varphi})} = \frac{1}{1-x-2x\cos{\varphi}}
\]
sein wird; und so werden wir sogar für $P$ gleich diesen endlichen Ausdruck haben:
\[
	P = \frac{1}{\pi}\int \frac{\partial \varphi}{1-x-2x\cos{\varphi}} \qquad \begin{bmatrix}
	\text{von}~ \varphi = 0 \\
	\text{bis}~ \varphi = \pi
	\end{bmatrix}
\]
welche Gleichung so beschafft werden kann:
\[
	P = \frac{1}{\pi (1-x)}\int \frac{\partial \varphi}{1 - \frac{2 x}{1-x}\cos{\varphi}} \qquad \begin{bmatrix}
	\text{von}~ \varphi = 0 \\
	\text{bis}~ \varphi = \pi
	\end{bmatrix}
\]
wo wir der Kürze wegen $\frac{2 x}{1-x} = k$ setzen wollen, sodass wir
\[
	P = \frac{1}{\pi (1-x)}\int\frac{\partial \varphi}{1-k\cos{\varphi}}
\]
haben.
\paragraph{§77}
Es ist aber bekannt, dass das Integral dieser Formel $\frac{\partial \varphi}{1+n\cos{\varphi}}$ gleich
\[
	\frac{1}{\sqrt{1-nn}}\arccos{\frac{\cos{\varphi} + n}{1+n\cos{\varphi}}}
\]
ist; daher, wenn wir $-k$ anstelle von $n$ schreiben, erhalten wir für unseren Fall
\[
	P = \frac{1}{\pi (1-x)\sqrt{1-kk}}\arccos{\frac{\cos{\varphi}-k}{1-k\cos{\varphi}}},
\]
wo die Addition einer Konstante nicht nötig ist, weil dieser Ausdruck im Fall $\varphi = 0$ von selbst verschwindet. Wir wollen daher für die andere Grenze $\varphi = \pi$ setzen, woher
\[
	\cos{\varphi} = -1 \quad \text{und} \quad \arccos{\frac{cos{\varphi}-k}{1-k\cos{\varphi}}} = \arccos{(-1)} = \pi
\]
wird; wir werden also
\[
	P = \frac{1}{(1-x)\sqrt{1-kk}}
\]
haben, welcher Ausdruck, wegen $\frac{2x}{1-x}$, übergeht in diesen:
\[
	P = \frac{1}{\sqrt{1-2x-3xx}},
\]
genauso wie vorher.
\paragraph{§78}
Weil
\[
	1 - 2x - 3xx = (1-x)^2 - 4xx = (1+x)(1-3x)
\]
ist, folgt, dass unsere zu summierende Reihe in zwei Fällen unendlich groß wird; natürlich in dem einen Fall, in dem $x=-1$ ist, in dem anderen aber, in dem $x=\frac{1}{3}$ ist. Dann aber wird unsere Reihe eine endliche Summe haben, wann immer $x$ zwischen diesen Grenzen enthalten ist: $-1$ und $\frac{1}{3}$. Wenn aber $x$ außerhalb dieser Grenzen angenommen wird, dann wird die Summe immer imaginär sein. So wird man für $x = \frac{1}{4}$ diese Summation haben:
\[
	1 + \frac{1}{4} + \frac{3}{4^2} + \frac{7}{4^3} + \frac{19}{4^4} + \frac{51}{4^5} + \frac{141}{4^6} + \mathrm{etc.} = \frac{4}{\sqrt{5}}
\]
\section*{\boldmath Untersuchung der Summe der übrigen Reihen $Q,~ R,~ S,~ \mathrm{etc,}$ die in §6 erörtert worden sind \unboldmath}
\paragraph{§79}
Wir wollen von der Reihe $Q$, die
\[
	Q = xx + 2x^3 + 6x^4 + \cdots + qx^{n+1} + q'x^{n+2} + q''x^{n+3} + \mathrm{etc.}
\]
ist, aus beginnen, deren erster Term, $xx$, aus der Potenz $n=1$ entsteht; dort, wenn wir genauso die Reihe von $n=0$ aus beginnen wollen, der Term $0x$ vorangestellt werden muss. Für diese Reihe aber haben wir oben gezeigt, dass $q = \frac{1}{2}(p'-p)$ ist, woher die Summe dieser Reihe aus der ersten Reihe $p$ auf folgende Weise gefunden werden kann.
\paragraph{§80}
Weil
\[
	P = 1 + x + 3xx + \cdots + px^n + p'x^{n+1} + \mathrm{etc.}
\]
ist, wird
\[
	Px = x + xx + \cdots + px^{n+1} + \mathrm{etc.}
\]
sein, welche letzte Reihe von der ersten abgezogen
\[
	P(1-x) = 1 + 2xx + \cdots + (p' - p)x^{n+1} + \mathrm{etc.}
\]
übriglässt. Weil daher $p' - p= 2q$ ist, wird
\[
	P(1-x) = 1 + 2Q
\]
sein; und so wird die Summe dieser Reihe bekannt, weil
\[
	Q = \frac{P(1-x)-1}{2}
\]
ist. Gerade vorher haben wir aber gesehen, dass
\[
	P = \frac{1}{\sqrt{1-2x-3xx}}
\]
ist und so werden wir
\[
	Q = \frac{1-x-\sqrt{1-2x-3xx}}{2\sqrt{1-2x-3xx}}
\]
haben.
\paragraph{§81}
Wir wollen zur Reihe $R$ weitergehen, die sich so verhielt:
\[
	R = x^4 + 3x^5 + 10x^6 + \cdots + rx^{n+2} + r'x^{n+3} + \mathrm{etc},
\]
deren erster Term, $x^4$, aus der Potenz $n=2$ entsprang, woher die beiden vorangestellten Terme $0x^2 + 0x^3$ aufzunehmen sind; um deren Summe zu finden bemerke man, dass $r = q'-q-p$ ist. Daher gilt, wenn die folgenden Operationen ausgeführt werden:
\begin{alignat*}{6}
	  Q  &= &xx  & +2x^3 &+\cdots  &+qx^{n+1}  &+q'x^{n+2}  &+\mathrm{etc.} \\
	-Qx &= &    & -x^3   &-\cdots &- \cdots    &-qx^{n+2}   &-\mathrm{etc.} \\
  -Px^2  &=  -&xx & -x^3  & -\cdots &- \cdots   &-px^{n+2}   &-\mathrm{etc.}
\end{alignat*}
daher wird zusammengenommen
\[
	Q(1-x) - Pxx = x^4 + 3x^5 + \cdots + (q' - q - p)x^{n+2} + \mathrm{etc.} = R
\]
werden.
\paragraph{§82}
Auf diese Weise also haben wir die Summe $R$ durch die beiden vorhergehenden Reihen $Q$ und $P$ bestimmt; weil diese schon bekannt sind, haben wir auch die Summe der Reihe $R$ algebraisch durch eine bestimmte Funktion von $x$ ausgedrückt erhalten; wie diese angenehm entwickelt werden kann, werden wir anschließend zeigen.
\paragraph{§85}
Für die Reihe $S$, die so vorgelegt worden war:
\[
	S = x^6 + 4x^7 + 15x^8 + \cdots + sx^{n+3} + s'x^{n+4} + \mathrm{etc.}
\]
sind ihr drei verschwindende Terme angefügt anzusehen: natürlich $0x^3 + 0x^4 + 0x^5$, wenn wir freilich von der Potenz $n=0$ aus anfangen wollen. Oben haben wir aber gefunden, dass $s = r' - r - q$ ist, woher wir die folgenden Operationen anstellen wollen:
\begin{alignat*}{8}
	R & = & x^4 & +3x^5 & +10x^6 & +\cdots & +rx^{n+2} & +r'x^{n+3} & +\mathrm{etc.} \\
  -Rx & = &     & -x^5  & -3x^6  & -\cdots & -\cdots   & -rx^{n+3}  & -\mathrm{etc.} \\
 -Qxx & = -&x^4 & -2x^5 & -6x^6  & -\cdots & -\cdots   & -qx^{n+3}  & -\mathrm{etc.}
\end{alignat*}
nach Zusammenfassen welcher drei Reihen diese Reihe entsteht
\[
	x^6 + \cdots + sx^{n+3},
\]
welche die Reihe $S$ selbst ist. Deshalb wird die Summe dieser Reihe durch die beiden vorhergehenden Reihen $Q$ und $R$ so bestimmt, dass
\[
	S = R(1-x) - Qxx
\]
ist, deren Entwicklung sogar hinreichend leicht gefunden werden kann, wie bald gezeigt werden wird.
\paragraph{§84}
Auf dieselbe Weise wird die Reihe $T$ durch die beiden vorhergehenden Reihen $R$ und $S$ bestimmt werden, uns zwar auf diese Weise:
\begin{alignat*}{8}
	S &= &x^6 & +4x^7 & +15x^8 & +\cdots & +sx^{n+3} & +s'x^{n+4} & +\mathrm{etc.} \\
	-Sx &= & & -x^7 & -4x^8 & -\cdots & -\cdots & -sx^{n+4} & +\mathrm{etc.} \\
	-Rxx &= -&x^6 & -3x^7 & -10x^8 & -\cdots & -\cdots & -rx^{n+4} & -\mathrm{etc.}
\end{alignat*}
Weil also $s'-s-r = t$ ist, werden diese drei Reihen zusammengefasst
\[
	S(1-x) - Rxx = x^8 + \cdots + tx^{n+4} + \mathrm{etc.}
\]
ergeben; weil dies die Reihe $T$ selbst ist, wird
\[
	T = S(1-x) - Rxx
\]
sein.
\paragraph{§85}
Daher ist also klar, dass jede einzelne dieser Reihen hinreichend leicht durch die beiden vorhergehenden bestimmt werden kann und sogar durch ein völlig gleichmäßiges Bildungsgesetz. Wir wollen sie zusammengefasst vor die Augen stellen:
\begin{align*}
Q &= \frac{P(1-x)-1}{2} \\
R &= Q(1-x) - Pxx \\
S &= R(1-x) - Qxx \\
T &= S(1-x) - Rxx \\
U &= T(1-x) - Sxx
\end{align*}
woher klar ist, dass all diese Summen nach einer rekurrenten Reihe fortschreiten, deren Relationsskala $(1-x),\, -xx$ ist. Es wird aber bald klar werden, dass diese Reihe sogar eine geometrische ist.
\paragraph{§86}
Um das zu zeigen, wollen wir, weil nach Ausführung der Entwicklung
\[
	\frac{Q}{P} = \frac{1-x-\sqrt{1-2x-3xx}}{2}
\]
ist, der Kürze wegen
\[
	\frac{1-x-\sqrt{1-2x-3xx}}{2} = \nu
\]
setzen, sodass wir $Q = P\nu$ haben; daher aber wird, nach Wegschaffen der Irrationalität, weil
\[
	\sqrt{1-2x-3xx} = 1-x-2\nu
\]
ist, diese Gleichung entstehen:
\[
	(1-x)^2 - 4xx = (1-x)^2 - 4\nu (1-x) + 4\nu\nu,
\]
die zurückgeführt wird auf diese:
\[
	\nu (1-x) - xx = \nu\nu,
\]
was natürlich förderlich sein wird, es bemerkt zu haben.
\paragraph{§87}
Gleich wird für diese Reihe $R$, wenn wir anstelle von $Q$ diesen Wert $P\nu$ einsetzen, diese Gleichung entstehen:
\[
	R = P(\nu(1-x)-xx)
\]
und daher durch die gerade bemerkte Relation
\[
	R = P\nu\nu.
\]
Wenn wir weiter anstelle von $Q$ und $R$ die gefundenen Werte schreiben, werden wir auf die gleiche Weise
\begin{align*}
	S &= P\nu (\nu (1-x)-xx) = P\nu^3 \\
	T &= P\nu\nu (\nu (1-x) - xx) = P\nu^4 \\
	U &= P\nu^3 (\nu (1-x) -xx) = P\nu^5 \\
	\vdots \\
	Z &= P\nu^{\lambda -1} (\nu (1-x) -xx) = P\nu^{\lambda +1}
\end{align*}
erhalten.
\paragraph{§88}
Wenn wir also gleich diese Bestimmungen zu bekannten Integralformeln, welche wir für die Buchstaben $p,\, q,\, r,\, \mathrm{etc}$ gefunden haben, transformieren, weil wir ja
\[
	Z = \frac{1}{\pi}\int \partial \varphi \cos{\lambda\varphi}(1+2\cos{\varphi})^n
\]
gefunden haben, wird, wenn wir dem Exponenten $n$ nacheinander die Werte $0,\, 1,\, 2,\, 3,\, 4,\, \mathrm{etc}$ zuteilen, weil die Reihe $Z$ anzusehen ist von der Potenz $x^{\lambda}$ zu beginnen, die Differentialformel $\partial\varphi \cos{\lambda}\varphi$ mit dieser geometrischen Reihe multipliziert werden müssen:
\[
	(1+2\cos{\varphi})^0 x^{\lambda} + (1+2\cos{\varphi})^1 x^{\lambda +1} + (1 + 2\cos{\varphi})^2 x^{\lambda +2} + \mathrm{etc},
\]
deren Summe
\[
	\frac{x^{\lambda}}{1-x-2x\cos{\varphi}}
\]
ist, nach Einführung welcher in die Rechnung also die gesuchte Summe so ausgedrückt werden wird:
\[
	Z = \frac{1}{\pi}\int \frac{x^{\lambda} \partial \varphi \cos{\lambda \varphi}}{1-x-2x\cos{\varphi}} \qquad \begin{bmatrix}
	\text{von}~ \varphi = 0 \\
	\text{bis}~ \varphi = \pi
	\end{bmatrix}
\]
wo die Größe $x$ konstant ist.	
\paragraph{§89}
Weil wir ja also hier diese Summe gefunden haben, natürlich
\[
	Z = P\nu^{\lambda} = \frac{\nu^{\lambda}}{\sqrt{1-2x-3xx}},
\]
während
\[
	\nu = \frac{1 - x - \sqrt{1-2x-3xx}}{2}
\]
wird, werden wir nun den Wert dieser Integralformeln sogar algebraisch ausdrücken können, weil wir ja nun wissen, dass
\[
	\frac{1}{\pi}\int \frac{x^{\lambda} \partial \varphi \cos{\lambda \varphi}}{1-x-2x\cos{\varphi}} = \frac{\nu^{\lambda}}{\sqrt{1-2x-3xx}},
\]
oder wir werden durch Multiplizieren mit $\frac{\pi}{x^{\lambda}}$
\[
	\int\frac{\partial\varphi \cos{\lambda \varphi}}{1-x-2x\cos{\varphi}} = \frac{\pi}{\sqrt{1-2x-3xx}}\ebinom{\nu}{x}^{\lambda}
\]
haben.
\paragraph{§90}
Weil ja die Integration größerer Aufmerksamkeit würdig erscheint, wollen wir diese in eine gefälligere Form bringen und, weil ja $x$ und $\nu$ hier als Konstanten betrachtet werden, wollen wir $\frac{\nu}{x} = b$ setzen, und wegen
\[
	\nu = 1 - x - \sqrt{1-2x-3xx}
\]
wird
\[
	2bx = 1 - x - \sqrt{1-2x-3xx}
\]
sein, welche Gleichung, nachdem die Irrationalität weggeschafft worden ist,
\[
	4bbxx - 4bx(1-x) + (1-x)^2 = (1-x)^2 - 4xx
\]
liefert, die zurückgeführt wird auf diese:
\[
	bbx - b  + bx = -x,
\]
woher die Größe $x$ selbst hinreichend angenehm bestimmt wird, weil $x = \frac{b}{bb+b+1}$ wird und daher
\[
	1-x = \frac{bb+1}{bb+b+1},
\]
und daher wird weiter, weil
\[
	\sqrt{1-2x-3xx} = 1 - x - 2bx
\]
war, nun
\[
	\sqrt{1-2x-3xx} = \frac{1-bb}{1+b+bb}
\]
sein.
\paragraph{§91}
Wenn wir daher also anstelle der Größe $x$ den Buchstaben $b$ in unsere Rechnung einführen, wird die gefundene Integration auf diese einfachere Form zurückgeführt werden:
\[
	\int\frac{\partial \varphi \cos{\lambda \varphi}}{1-2b\cos{\varphi} +bb} \qquad \begin{bmatrix}
	\text{von}~ \varphi = 0 \\
	\text{bis}~ \varphi = \pi
	\end{bmatrix}
	 = \frac{\pi b^{\lambda}}{1-bb}
\]
deren Gültigkeit aus den bisher gefundenen Rechnungen gefolgert worden ist; sie kann aber auch sofort und direkt gezeigt werden, wodurch alles Vorhergehende umso mehr bestätigt werden wird.
\paragraph{§92}
Um das also zu zeigen, wollen wir die hinreichend bekannte Reduktion zur Hilfe nehmen, durch die
\[
	\int \frac{\partial \varphi}{\alpha + \beta\cos{\varphi}} = \frac{1}{\sqrt{\alpha\alpha - \beta\beta}}\arccos{\frac{\alpha \cos{\varphi} + \beta}{\alpha + \beta\cos{\varphi}}}
\]
ist. Es werde nun $\alpha = 1 + bb$ und $\beta = -2b$ und wir werden
\[
	\int \frac{\partial \varphi}{1-2b\cos{\varphi} + bb} = \frac{1}{1-bb}\arccos{\frac{(1+bb)\cos{\varphi} - 2b}{1-2b\cos{\varphi} + bb}}
\]
haben, welches Integral gleich für $\varphi = 0 $ gesetzt verschwindet. Nachdem also für die andere Grenze $\varphi = \pi$ gesetzt wurde, wird diese Integral $\frac{\pi}{1-bb}$ werden.
\paragraph{§93}
Weil wir ja also für unsere Integrationsgrenzen
\[
	\int\frac{\partial\varphi}{1-2b\cos{\varphi} + bb} = \frac{\pi}{1-bb}
\]
gefunden haben und natürlich
\[
	\int\partial\varphi = \pi \quad \text{und daher} \quad \int\frac{\partial\varphi (1-2b\cos{\varphi} +bb)}{1-2b\cos{\varphi} + bb} = \pi
\]
ist, werden wir, indem wir diese Form in zwei Teile aufteilen,
\[
	\pi = (1+bb)\int\frac{\partial \varphi}{1-2b\cos{\varphi}+bb} - 2b\int\frac{\partial \varphi \cos{\varphi}}{1-2b\cos{\varphi} + bb}
\]
haben, woher wir
\[
	\int\frac{\partial\varphi\cos{\varphi}}{1-2b\cos{\varphi}+bb} = \frac{\pi b}{1-bb}
\]
berechnen.
\paragraph{§94}
Weil ja für unsere Integrationsgrenzen im allgemeinen
\[
	\int\partial\varphi\cos{i\varphi} = 0
\]
ist, wenn natürlich $i$ eine ganze Zahl war, wollen wir diese Formel mit $1+bb-2b\cos{\varphi}$ erweitern und werden
\[
	\int\frac{\partial\varphi((1+bb)\cos{i\varphi} - b\cos{(i-1)\varphi} - b\cos{(i+1)\varphi})}{1+bb-2b\cos{\varphi}}
\]
erhalten. Diese Form wird uns in drei Teile geteilt
\[
	(1+bb)\int\frac{\partial \varphi \cos{i\varphi}}{1-2b\cos{\varphi} + bb} = b\int\frac{\partial\varphi\cos{(i-1)\varphi}}{1-2b\cos{\varphi} + bb} + b\int\frac{\partial \varphi\cos{(i+1)\varphi}}{1-2b\cos{\varphi}+bb}
\]
geben, woher wir diese allgemeine Reduktion berechnen:
\[
	\int\frac{\partial \varphi\cos{(i+1)\varphi}}{1-2b\cos{\varphi}+bb} = \frac{1+bb}{b}\int\frac{\partial \varphi \cos{i\varphi}}{1-2b\cos{\varphi}+bb} - \int\frac{\partial\varphi \cos{(i-1)\varphi}}{1-2b\cos{\varphi}+bb}
\]
mit deren Hilfe aus den beiden Integrale für die Winkel $i\varphi$ und $(i-1)\varphi$ das Integral für den Winkel $(i+1)\varphi$ bestimmt werden kann, woher sich die folgende Tabelle aufstellen lassen wird:
\begin{align*}
\int\frac{\partial\varphi}{1+bb-2b\cos{\varphi}} &= \frac{\pi}{1-bb} \\
\int\frac{\partial\varphi\cos{\varphi}}{1+bb-2b\cos{\varphi}} &= \frac{\pi b}{1-bb} \\
\int\frac{\partial\varphi\cos{2\varphi}}{1+bb-2b\cos{\varphi}} &= \frac{\pi bb}{1-bb} \\
\int\frac{\partial\varphi\cos{3\varphi}}{1+bb-2b\cos{\varphi}} &= \frac{\pi b^3}{1-bb} \\
\int\frac{\partial\varphi\cos{4\varphi}}{1+bb-2b\cos{\varphi}} &= \frac{\pi b^4}{1-bb} \\
\vdots \\
\int\frac{\partial\varphi\cos{\lambda\varphi}}{1+bb-2b\cos{\varphi}} &= \frac{\pi b^{\lambda}}{1-bb}
\end{align*}
völlig wie wir es oben gefunden haben.
\end{document}